\documentclass[a4paper,11pt]{amsart}

\usepackage{graphicx}
\usepackage{mathptmx}
\usepackage{amsmath}
\usepackage{amssymb}
\usepackage{enumitem}
\usepackage{xcolor}

\newmuskip\pFqmuskip

\newcommand*\pFq[6][8]{%
  \begingroup 
  \pFqmuskip=#1mu\relax
  \mathcode`=\string"8000
  \begingroup\lccode`\~=`\,
  \lowercase{\endgroup\let~}\pFqcomma
  F^{#2}_{#3}{\left(\genfrac..{0pt}{}{#4}{#5}\bigg|#6\right)}%
  \endgroup
}
\newcommand{\pFqcomma}{\mskip\pFqmuskip}

\newtheorem{theorem}{Theorem}[section]

\newtheorem{corollary}[theorem]{Corollary}
\newtheorem{proposition}[theorem]{Proposition}

\newtheorem{remark}[theorem]{Remark}

\begin{document}

\title[Probabilistic Stirling and degenerate Stirling numbers]{Probabilistic Stirling numbers associated with sequences }

\author{Dae San  Kim }
\address{Department of Mathematics, Sogang University, Seoul 121-742, Republic of Korea}
\email{dskim@sogang.ac.kr}

\author{Taekyun  Kim}
\address{Department of Mathematics, Kwangwoon University, Seoul 139-701, Republic of Korea}
\email{tkkim@kw.ac.kr}


\begin{abstract}
Let $Y$ be a random variable whose moment generating function exists in a neighborhood of the origin. Recently, probabilistic Stirling numbers of the first kind, $s_{Y}(n,k)$, and of the second kind, $S_{2}^{Y}(n,k)$, associated with $Y$ have been introduced. However, $s_{Y}(n,k)$, based on the cumulant generating function of $Y$, and $S_{2}^{Y}(n,k)$ do not satisfy orthogonality and inverse relations. This paper aims to redefine the probabilistic Stirling numbers of the first kind associated with $Y$, denoting them by $S_{1}^{Y}(n,k)$, such that $S_{1}^{Y}(n,k)$ and $S_{2}^{Y}(n,k)$ do satisfy these crucial relations. Furthermore, we investigate their degenerate counterparts, the probabilistic degenerate Stirling numbers of both kinds. We explicitly compute $S_{2}^{Y}(n,k), S_{1}^{Y}(n,k)$, and their degenerate versions for several discrete and continuous random variables $Y$. As an application, we express arbitrary polynomials as linear combinations of probabilistic Euler and probabilistic degenerate Euler polynomials.

\end{abstract}

\maketitle

\section{\bf Introduction}
The Stirling number of the second, denoted by $S_{2}(n, k)$, counts the number of ways to partition a set of $n$ objects into $k$ nonempty subsets. The signed Stirling number of the first kind, denoted by, $S_{1}(n, k)$ is defined such that the number of permutations of $n$ elements with exactly $k$ cycles is given by $|S_{1}(n,k)|=(-1)^{n-k}S_{1}(n,k)$ (see, for example, [7,23,24,26]).
The degenerate Stirling numbers of the second kind, $S_{2,\lambda}(n, k)$ (defined in (3)), and of the first kind, $S_{1,\lambda}(n, k)$ (defined in (5)), are degenerate versions of $S_{2}(n,k)$ and $S_{1}(n,k)$, respectively (see [13,14,22]). These degenerate Stirling numbers frequently arise in the study of degenerate versions of various special polynomials and numbers (see [8,9,20]). \par
The study of degenerate versions of special polynomials and numbers has seen a resurgence of interest and a corresponding increase in research activity. These investigations have employed a diverse range of mathematical tools, including quantum mechanics, generating functions, combinatorial methods, $p$-adic analysis, umbral calculus, probability theory, differential equations, analytic number theory, and operator theory. This exploration of degenerate versions originated with Carlitz's work on degenerate Bernoulli and Euler polynomials (see [4,5]). Notably, this area of study has expanded beyond polynomials and numbers to encompass transcendental functions (see [15,16]) and umbral calculus (see [10]). \par
The Stirling numbers of both kinds and the degenerate Stirling numbers of both kinds are generalized as follows (see [11]). Let ${\bf{P}}=\{p_{n}(x)\}_{n=0}^{\infty}$ be a sequence of polynomials such that deg\,$p_{n}(x)=n,\, p_{0}(x)=1$. The Stirling numbers of the second kind associated with $\bf{P}$, denoted by $S_2(n,k;{\bf{P}})$, are defined as the coefficients in the expansion of $p_{n}(x)$ in terms of the falling factorials $(x)_{k}$:
$p_{n}(x)=\sum_{k=0}^{n}S_2(n,k;{\bf{P}})(x)_{k}$.
Note that $S_{2}(n,k;{\bf{P}})$ reduces to the ordinary Stirling numbers of the second kind, $S_{2}(n,k)$, when ${\bf{P}}=\{x^n\}_{n=0}^{\infty}$ (see \eqref{1}), and to the degenerate Stirling numbers of the second kind, $S_{2,\lambda}(n, k)$, when ${\bf{P}}=\{(x)_{n,\lambda}\}_{n=0}^{\infty}$ (see \eqref{3}).
The Stirling numbers of the first kind associated with $\bf{P}$, denoted by $S_{1}(n,k;{\bf{P}})$, are defined by the expansion $(x)_{n}=\sum_{k=0}^{n}S_{1}(n,k;{\bf{P}})p_{k}(x)$. The Stirling numbers of both kinds associated with $\bf{P}$ have been investigated using umbral calculus (see [11]). Analogously, the central factorial numbers of both kinds associated with $\bf{P}$ can be studied by replacing the falling factorials $(x)_{n}$ with the central factorials $x^{[n]}$, where
$x^{[n]} = x(x+\frac{1}{2}n-1)_{n-1},\,\, (n \ge 1),\,\, x^{[0]}=1$,\,\, (see [12]).  \par
Let $Y$ be a random variable whose moment generating function exists in some neighborhood of the origin (see \eqref{8}). Recently, Adell and Lekuona introduced the probabilistic Stirling numbers of the second kind associated with $Y$, denoted by $S_{2}^{Y}(n,k)$ (see [1,3]). Furthermore, Adell and B\'{e}nyi defined the probabilistic Stirling numbers of the first kind associated with $Y$, denoted by $s_{Y}(n,k)$, using the cumulant generating function (see [2]). However, $S_{2}^{Y}(n,k)$ and $s_{Y}(n,k)$ do not satisfy orthogonality and inverse relations. In addition, $s_{Y}(n,k)$ does not reduce to the classical Stirling numbers of the first kind, $S_{1}(n,k)$, when $Y=1$. We address these shortcomings by redefining the probabilistic Stirling numbers of the first kind, denoting them by $S_{1}^{Y}(n,k)$. Our redefined $S_{1}^{Y}(n,k)$ and $S_{2}^{Y}(n,k)$ do satisfy orthogonality and inverse relations (see Proposition 1.1), and $S_{1}^{Y}(n,k)$ reduces to $S_{1}(n,k)$, when $Y=1$.
We also consider the probabilistic degenerate Stirling numbers of the second kind associated with $Y$, denoted by $S_{2,\lambda}^{Y}(n,k)$ (see \eqref{14} and [21]), and the probabilistic degenerate Stirling numbers of the first kind associated with $Y$, denoted by $S_{1,\lambda}^{Y}(n,k)$ (see \eqref{15}). $S_{2,\lambda}^{Y}(n,k)$ and $S_{1,\lambda}^{Y}(n,k)$ satisfy orthogonality and inverse relations (see Proposition 1.2) and $S_{1,\lambda}^{Y}(n,k)$ reduces to the standard degenerate Stirling numbers of the first kind, $S_{1,\lambda}(n,k)$, when $Y=1$.
However, in [18], $S_{1,\lambda}^{Y}(n,k)$ is defined differently in that work using the degenerate cumulant generating function. Consequently, the $S_{1,\lambda}^{Y}(n,k)$ defined there and $S_{2,\lambda}^{Y}(n,k)$ do not together satisfy orthogonality and inverse relations. In Section 2, we explicitly compute $S_{2}^{Y}(n,k), S_{1}^{Y}(n,k), S_{2,\lambda}^{Y}(n,k)$, and $S_{1,\lambda}^{Y}(n,k)$ for several discrete and continuous random variables $Y$ (see Theorems 2.1-2.8).
Section 3 illustrates the orthogonality relations for the geometric random variable with parameter $0 <p <1$. Furthermore, we address the problem of expressing an arbitrary polynomial in terms of the probabilistic Euler polynomials $\mathcal{E}_{n}^{Y}(x)$ and the probabilistic degenerate Euler polynomials $\mathcal{E}_{n,\lambda}^{Y}(x)$, demonstrating the necessity of inversion relations. Probabilistic extensions of numerous special polynomials and numbers have been studied recently (see, for example, [6,17,19,28,29] and the references therein).

\vspace{0.1in}

The Stirling numbers of the second kind, $S_{2}(n,k)$, are defined by
\begin{align}
&x^{n}=\sum_{k=0}^{n}S_{2}(n,k)(x)_{k},\label{1} \\
&\frac{1}{k!}(e^{t}-1)^{k}=\sum_{n=k}^{\infty}S_{2}(n,k)\frac{t^{n}}{n!}, \nonumber\\
&S_{2}(n,k)=\frac{1}{k!}\sum_{j=0}^{k}\binom{k}{j}(-1)^{k-j}j^{n}, \nonumber
\end{align}
where the falling factorials $(x)_{n}$ are given by
\begin{equation*}
(x)_{0}=1, \quad (x)_{n}=x(x-1)\cdots(x-n+1),\quad (n \ge 1).
\end{equation*}
The Stirling numbers of the first kind, $S_{1}(n,k)$, are defined by
\begin{align}
&(x)_{n}=\sum_{k=0}^{n}S_{1}(n,k)x^{k}, \label{2} \\
&\frac{1}{k!}\big(\log(1+t)\big)^{k}=\sum_{n=k}^{\infty}S_{1}(n,k)\frac{t^{n}}{n!}. \nonumber
\end{align}
Let $\lambda$ be any nonzero real number. The degenerate Stirling numbers of the second kind, $S_{2,\lambda}(n,k)$, are defined by
\begin{align}
&(x)_{n,\lambda}=\sum_{k=0}^{n}S_{2,\lambda}(n,k)(x)_{k}, \label{3} \\
&\frac{1}{k!}(e_{\lambda}(t)-1)^{k}=\sum_{n=k}^{\infty}S_{2,\lambda}(n,k)\frac{t^{n}}{n!}, \nonumber \\
&S_{2,\lambda}(n,k)=\frac{1}{k!}\sum_{j=0}^{k}\binom{k}{j}(-1)^{k-j}(j)_{n,\lambda},\nonumber
\end{align}
where $(x)_{n,\lambda}$ are the degenerate (or generalized) falling factorials given by
\begin{equation*}
(x)_{0,\lambda}=1,\ (x)_{n,\lambda}=x(x-\lambda)\cdots(x-(n-1)\lambda),\ (n\ge 1),
\end{equation*}
and $e_{\lambda}^{x}(t)$ are the degenerate exponentials defined by
\begin{equation}
e_{\lambda}^{x}(t)=(1+\lambda t)^{\frac{x}{\lambda}}=\sum_{n=0}^{\infty}(x)_{n,\lambda}\frac{t^{n}}{n!}, \quad e_{\lambda}(t)=e_{\lambda}^{1}(t),\quad (\mathrm{see}\ [9,17,18,21]). \label{4}
\end{equation}
Note here that $\lim_{\lambda\rightarrow 0}e_{\lambda}^{x}(t)=e^{xt}$. \\
The degenerate Stirling numbers of the first kind, $S_{1,\lambda}(n,k)$, are given by
\begin{align}
&(x)_{n}=\sum_{k=0}^{n}S_{1,\lambda}(n,k)(x)_{k,\lambda}, \label{5} \\
&\frac{1}{k!}\big(\log_{\lambda}(1+t)\big)^{k}=\sum_{n=k}^{\infty}S_{1,\lambda}(n,k)\frac{t^{n}}{n!}, \nonumber
\end{align}
where $\log_{\lambda}(t)$ are the degenerate logarithm defined by
\begin{equation}
\log_{\lambda}(t)=\frac{1}{\lambda}\big(t^{\lambda}-1\big). \label{6}
\end{equation}
Note here that the degenerate exponential $e_{\lambda}(t)$ in \eqref{4} and the degenerate logarithm $\log_{\lambda}(t)$ in \eqref{6} satisfy
\begin{equation}
e_{\lambda}\big(\log_{\lambda}(t) \big)=\log_{\lambda}\big(e_{\lambda}(t) \big)=t. \label{7}
\end{equation}
Among other things, the Stirling numbers of both kinds and the degenerate Stirling numbers of both kinds satisfy orthogonality and inversion relations (see Propositions 1.1 and 1.2 with $Y=1$). \par
Let $Y$ be a random variable whose moment generating function exists in a neighborhood of the origin:
\begin{equation}
E[e^{Yt}]=\sum_{n=0}^{\infty}E[Y^{n}]\frac{t^{n}}{n!}\,\,\, \mathrm{exists,\,\, for}\,\,|x|<r, \label{8}
\end{equation}
for some positive real number $r$. \\
Let $(Y_{j})_{j \ge 1}$ be a sequence of mutually independent copies of the random variable $Y$, and let
\begin{equation}
S_{k}=Y_{1}+Y_{2}+\cdots +Y_{k}, \quad (k \ge 1), \quad  S_{0}=0. \label{9}
\end{equation}
The probabilistic Stirling numbers of the second kind associated with $Y$, $S_{2}^{Y}(n,k)$, are given by (see \eqref{9})
\begin{align}
&\frac{1}{k!}(E[e^{Yt}]-1)^{k}=\sum_{n=k}^{\infty}S_{2}^{Y}(n,k)\frac{t^{n}}{n!}, \label{10}\\
&S_{2}^{Y}(n,k)=\frac{1}{k!}\sum_{j=0}^{k}\binom{k}{j}(-1)^{k-j}E[S_{j}^{n}]. \nonumber
\end{align}
From the definition in \eqref{10}, it is immediate to see that
\begin{equation}
S_{2}^{Y}(k,k)=E[Y]^{k}. \label{11}
\end{equation}
Here we define the probabilistic Stirling numbers of the first kind associated with $Y$, $S_{1}^{Y}(n,k)$, in a different way from those ones in [2] and [18]. We believe that they should be defined in such a manner that they obey orthogonality and inverse relations. To do this, we assume from now on till the end of this paper that
\begin{equation}
E[Y] \ne 0. \label{12}
\end{equation}
We introduce the notation:
\begin{equation}
e_{Y}(t)=E[e^{Yt}]-1. \label{13}
\end{equation}
Then we have
\begin{equation}
\frac{1}{k!}\big(e_{Y}(t)\big)^{k}=\sum_{n=k}^{\infty}S_{2}^{Y}(n,k)\frac{t^{n}}{n!}. \label{14}
\end{equation}
If $f(t)=\sum_{n=0}^{\infty}a_{n}\frac{t^{n}}{n!}$ is a delta series, namely $a_{0}=0$ and $a_{1} \ne 0$, then the compositional inverse $\bar{f}(t)$ of $f(t)$ satisfying $f(\bar{f}(t))=\bar{f}(f(t))=t$ exists. For example, $e_{\lambda}(t)$ and $\log_{\lambda}(t)$ are compositional inverses to each other (see \eqref{7}). Note that, as $e_{Y}(t)=E[Y]t+\sum_{m=2}^{\infty}E[Y^{m}]\frac{t^{m}}{m!}$ and $E[Y] \ne 0$ (see \eqref{12}, \eqref{13}), $e_{Y}(t)$ is a delta series. \par
Now, we define the  probabilistic Stirling numbers of the first kind associated with $Y$ by: for $k \ge 0$,
\begin{equation}
\frac{1}{k!}\big(\bar{e}_{Y}(t)\big)^{k}=\sum_{n=k}^{\infty}S_{1}^{Y}(n,k)\frac{t^{n}}{n!}, \label{15}
\end{equation}
where $\bar{e}_{Y}(t)$ is the compositional inverse of $e_{Y}(t)$. \\
In addition, as usual, we agree that
\begin{equation}
S_{2}^{Y}(n,k)=S_{1}^{Y}(n,k)=0,\,\, \mathrm{if}\,\, k>n\,\,\mathrm{or}\,\, k <0. \label{16}
\end{equation}
Note that $S_{2}^{Y}(n,k)=S_{2}(n,k),\,\, S_{1}^{Y}(n,k)=S_{1}(n,k)$, when $Y=1$. \par
Using the definitions in \eqref{14} and \eqref{15}, one shows that $S_{2}^{Y}(n,k)$ and $S_{1}^{Y}(n,k)$ satisfy the orthogonality relations in (a) of Proposition 1.1, from which the inverse relations in (b) and (c) follow.
\begin{proposition}
The following orthogonality and inverse relations are valid for $S_{1}^{Y}(n,k)$ and $S_{2}^{Y}(n,k)$.
\begin{flalign*}
&(a)\,\, \,\sum_{k=l}^{n} S_{2}^{Y}(n,k)S_{1}^{Y}(k,l)=\delta_{n,l}, \quad \sum_{k=l}^{n} S_{1}^{Y}(n,k)S_{2}^{Y}(k,l)=\delta_{n,l}, \\
&(b)\,\, a_{n}=\sum_{k=0}^{n}S_{2}^{Y}(n,k) b_{k}\,\, \iff \,\, b_{n}=\sum_{k=0}^{n}S_{1}^{Y}(n,k)a_{k}, \\
&(c)\,\, a_{n}=\sum_{k=n}^{m}S_{2}^{Y}(k,n)b_{k} \,\, \iff \,\, b_{n}=\sum_{k=n}^{m}S_{1}^{Y}(k,n)a_{k}. &&
\end{flalign*}
\end{proposition}
Now, we introduce the following:
\begin{equation}
\bar{f}_{Y}(t)=\log E[e^{Yt}]. \label{17}
\end{equation}
Let $f_{Y}(t)$ be the compositional inverse of $\bar{f}_{Y}(t)$. Here we have to observe that $\bar{f}_{Y}(t)$ is a delta series (see \eqref{18}). Indeed, we have (see \eqref{10}-\eqref{12})
\begin{align}
\bar{f}_{Y}(t)&=\log \big(1+(E[e^{Yt}]-1)\big) \label{18}\\
&=\sum_{j=1}^{\infty}(-1)^{j-1}(j-1)! \frac{1}{j!}(E[e^{Yt}]-1)^{j} \nonumber \\
&=\sum_{j=1}^{\infty}(-1)^{j-1}(j-1)!\sum_{n=j}^{\infty}S_{2}^{Y}(n,j)\frac{t^{n}}{n!} \nonumber \\
&=S_{2}^{Y}(1,1)t+\sum_{n=2}^{\infty}\sum_{j=1}^{n}(-1)^{j-1}(j-1)!S_{2}^{Y}(n,j)\frac{t^{n}}{n!} \nonumber \\
&=E[Y]t+\sum_{n=2}^{\infty}\sum_{j=1}^{n}(-1)^{j-1}(j-1)!S_{2}^{Y}(n,j)\frac{t^{n}}{n!}. \nonumber
\end{align}
Here we recall that
\begin{equation}
\bar{f}_{Y}(t)=\log E[e^{Yt}]=\sum_{n=1}^{\infty}\kappa_{n}(Y)\frac{t^{n}}{n!} \label{19}
\end{equation}
is called the cumulant generating function of $Y$. From \eqref{18} and \eqref{19}, we derive
\begin{equation*}
\kappa_{n}(Y)=\sum_{j=1}^{n}(-1)^{j-1}(j-1)!S_{2}^{Y}(n,j), \quad (n \ge 1).
\end{equation*}
In [2], the probabilistic Stirling numbers of the first kind associated with $Y$, $s_{Y}(n,k)$, are defined as
\begin{equation}
\frac{1}{k!}\big(\log E[e^{Yt}] \big)^{k}=\sum_{n=k}^{\infty}(-1)^{n-k}s_{Y}(n,k)\frac{t^{n}}{n!}. \label{20}
\end{equation}
However, $s_{Y}(n,k)$ and $S_{2}^{Y}(n,k)$ do not together satisfy orthogonality and inverse relations, and $s_{Y}(n,k)=\delta_{n,k}$,\,\, when $Y=1$. Here we note, in passing, that
\begin{equation*}
\kappa_{1}(Y)=E[Y],\,\,\kappa_{2}(Y)=\mathrm{Var}(Y),\,\, \kappa_{3}(Y)=E[(Y-E(Y))^{3}].
\end{equation*} \par
The probabilistic degenerate Stirling numbers of the second kind associated $Y$, $S_{2, \lambda}^{Y}(n,k)$, are defined by
\begin{align}
&\frac{1}{k!}\big(E[e_{\lambda}^{Y}(t)]-1\big)^{k}=\sum_{n=k}^{\infty}S_{2,\lambda}^{Y}(n,k)\frac{t^{n}}{n!}, \label{21} \\
&S_{2,\lambda}^{Y}(n,k)=\frac{1}{k!}\sum_{j=0}^{k}\binom{k}{j}(-1)^{k-j}E[(S_{j})_{n,\lambda}]. \nonumber
\end{align}
To define the probabilistic degenerate Stirling numbers of the first kind associated with $Y$, $S_{1,\lambda}^{Y}(n,k)$, we let
\begin{equation}
e_{Y,\lambda}(t)=E[e_{\lambda}^{Y}(t)]-1. \label{22}
\end{equation}
Then we have
\begin{equation}
\frac{1}{k!}(e_{Y,\lambda}(t))^{k}=\sum_{n=k}^{\infty}S_{2,\lambda}^{Y}(n,k)\frac{t^{n}}{n!}. \label{23}
\end{equation}
Note that $e_{Y,\lambda}(t)=E[Y]t+\sum_{m=2}^{\infty}E[(Y)_{m,\lambda}]\frac{t^{m}}{m!}$ is a delta series (see \eqref{12}). Now, we define the probabilistic degenerate Stirling numbers of the first kind associated with $Y$, $S_{1,\lambda}^{Y}(n,k)$, by
\begin{equation}
\frac{1}{k!}\big(\bar{e}_{Y,\lambda}(t)\big)^{k}=\sum_{n=k}^{\infty}S_{1,\lambda}^{Y}(n,k)\frac{t^{n}}{n!}, \label{24}
\end{equation}
where $\bar{e}_{Y,\lambda}(t)$ is the compositional inverse of $e_{Y,\lambda}(t)$.
Then $S_{2,\lambda}^{Y}(n,k)$ and $S_{1,\lambda}^{Y}(n,k)$ satisfy orthogonality and inverse relations, and $S_{2,\lambda}^{Y}(n,k)=S_{2,\lambda}(n,k)$ and $S_{1,\lambda}^{Y}(n,k)=S_{1,\lambda}(n,k)$, when $Y=1$. We also consider
\begin{equation}
\bar{f}_{Y,\lambda}(t)=\log E[e_{\lambda}^{Y}(t)],  \label{25}
\end{equation}
and its compositional inverse $f_{Y,\lambda}(t)$ as well.
As for the probabilisitc Stirling numbers associated with $Y$, the orthogonality and inverse relations are valid for the probabilistic degenerate Stirling numbers associated with $Y$.
\begin{proposition}
The following orthogonality and inverse relations are valid for $S_{1,\lambda}^{Y}(n,k)$ and $S_{2,\lambda}^{Y}(n,k)$.
\begin{flalign*}
&(a)\,\, \sum_{k=l}^{n} S_{2,\lambda}^{Y}(n,k)S_{1,\lambda}^{Y}(k,l)=\delta_{n,l}, \quad \sum_{k=l}^{n} S_{1,\lambda}^{Y}(n,k)S_{2,\lambda}^{Y}(k,l)=\delta_{n,l},\\
&(b)\,\, a_{n}=\sum_{k=0}^{n}S_{2,\lambda}^{Y}(n,k) b_{k}\,\, \iff \,\, b_{n}=\sum_{k=0}^{n}S_{1,\lambda}^{Y}(n,k)a_{k}, \\
&(c)\,\, a_{n}=\sum_{k=n}^{m}S_{2,\lambda}^{Y}(k,n)b_{k} \,\, \iff \,\, b_{n}=\sum_{k=n}^{m}S_{1,\lambda}^{Y}(k,n)a_{k}. &&
\end{flalign*}
\end{proposition}
The probabilistic Euler polynomials associated with $Y$, $\mathcal{E}_{n}^{Y}(x)$, are defined by
\begin{equation}
\frac{2}{E[e^{Yt}]+1}\big(E[e^{Yt}]\big)^{x}=\sum_{n=0}^{\infty}\mathcal{E}_{n}^{Y}(x)\frac{t^{n}}{n!}. \label{26}
\end{equation}
Moreover, the probabilistic degenerate Euler polynomials associated with $Y$, $\mathcal{E}_{n,\lambda}^{Y}(x)$, are given by
\begin{equation}
\frac{2}{E[e_{\lambda}^{Y}(t)]+1}\big(E[e_{\lambda}^{Y}(t)]\big)^{x}=\sum_{n=0}^{\infty}\mathcal{E}_{n,\lambda}^{Y}(x)\frac{t^{n}}{n!}. \label{27}
\end{equation} \par
Let $a(x)$ be a polynomial of degree of $n$ with complex coefficients. Then, for any integer $r$ with $0 \le r \le n$, (see [24])
\begin{align}
\Delta^{r}a(x)&=\sum_{i=0}^{r}\binom{r}{i}(-1)^{r-i}a(x+i) \label{28} \\
&=r!\sum_{j=r}^{n}S_{2}(j,r)\frac{1}{j!}a^{(j)}(x), \nonumber
\end{align}
where $\Delta$ is the forward difference operator given by (see [7,23,25])
\begin{equation}
\Delta a(x)=a(x+1)-a(x). \label{29}
\end{equation}
It is easy to see that
\begin{equation}
\Delta (x)_{n}=(x+1)_{n}-(x)_{n}=n(x)_{n-1},\quad (n \ge 0), \label{30}
\end{equation}
and hence, for any $r$ with $0 \le r \le n$,
\begin{equation}
\Delta^{r}(x)_{n}=(n)_{r}(x)_{n-r}. \label{31}
\end{equation}

\section{\bf Probabilistic Stirling and degenerate Stirling numbers}
In this section, for several discrete and continuous random variables $Y$, we determine
\begin{align*}
&e_{Y}(t), \,\,\, \bar{e}_{Y}(t), \,\,\, S_{2}^{Y}(n,k), \,\,\, S_{1}^{Y}(n,k),\,\,\, \bar{f}_{Y}(t), \,\,\, f_{Y}(t), \\
&e_{Y,\lambda}(t), \,\,\, \bar{e}_{Y,\lambda}(t), \,\,\, S_{2,\lambda}^{Y}(n,k), \,\,\, S_{1,\lambda}^{Y}(n,k),\,\,\, \bar{f}_{Y,\lambda}(t), \,\,\,\mathrm{and}\,\,\, f_{Y,\lambda}(t).
\end{align*}
Our reference here is [27], especially Table 2.1 on p. 61 and Table 2.2 on p. 62.\par
If $Y$ is a discrete random variable with probability mass function
$p(y)=P\{Y=y\}$, then for any real-valued function $g$,
\begin{equation}
E[g(Y)]=\sum_{y}g(y)p(y). \label{32}
\end{equation}
If $Y$ is a continuous random variable with probability density function $f (x)$, then
for any real-valued function $g$,
\begin{equation}
E[g(Y)]=\int_{-\infty}^{\infty}g(y)f(y) dy. \label{33}
\end{equation}
(a) Let $Y$ be the Bernoulli random variable. Then the probability mass function of $Y$ is given by
\begin{equation}
p(0)=1-p,\quad p(1)=p, \quad (0< p \le 1). \label{34}
\end{equation}
Then we have
\begin{align}
&E[Y]=p,\,\, E[e^{Yt}]=1+p(e^{t}-1),\,\,e_{Y}(t)=p(e^{t}-1),\,\,\bar{e}_{Y}(t)=\log(1+\frac{t}{p}), \label{35}\\
&\bar{f}_{Y}(t)=\log\big(1+p(e^{t}-1)\big),\,\,f_{Y}(t)=\log\big(1+\frac{1}{p}(e^{t}-1)\big). \nonumber
\end{align}
From \eqref{35}, we obtain
\begin{equation}
S_{2}^{Y}(n,k)=p^{k}S_{2}(n,k),\,\,S_{1}^{Y}(n,k)=\frac{1}{p^{n}}S_{1}(n,k). \label{36}
\end{equation}
We observe from \eqref{32} and \eqref{34} that
\begin{equation}
E[e_{\lambda}^{Y}(t)]=e_{\lambda}^{0}(t)(1-p)+e_{\lambda}(t)p=1+p(e_{\lambda}(t)-1). \label{37}
\end{equation}
From \eqref{37}, we derive
\begin{align}
&e_{Y,\lambda}(t)=p\big(e_{\lambda}(t)-1),\,\,\bar{e}_{Y,\lambda}(t)=\log_{\lambda}(1+\frac{t}{p}), \label{38}\\
&\bar{f}_{Y,\lambda}(t)=\log\big(1+p(e_{\lambda}(t)-1)\big),\,\,f_{Y,\lambda}(t)=\log_{\lambda}\big(1+\frac{1}{p}(e^{t}-1)\big). \nonumber
\end{align}
From \eqref{38}, we get
\begin{equation}
S_{2,\lambda}^{Y}(n,k)=p^{k}S_{2,\lambda}(n,k),\,\,S_{1,\lambda}^{Y}(n,k)=\frac{1}{p^{n}}S_{1,\lambda}(n,k). \label{39}
\end{equation}
We summarize the results in \eqref{36} ane \eqref{39} in the next theorem.
\begin{theorem}
Let $Y$ be the Bernoulli random variable with parameter $p$ with $0 < p \le 1$. Then, for\, $n \ge k$, we have
\begin{align*}
&S_{2}^{Y}(n,k)=p^{k}S_{2}(n,k),\,\,S_{1}^{Y}(n,k)=\frac{1}{p^{n}}S_{1}(n,k), \\
&S_{2,\lambda}^{Y}(n,k)=p^{k}S_{2,\lambda}(n,k),\,\,S_{1,\lambda}^{Y}(n,k)=\frac{1}{p^{n}}S_{1,\lambda}(n,k).
\end{align*}
\end{theorem}
\vspace{0.1in}
(b) Let $Y$ be the binomial random variable with parameter $(m,p)$. Then the probability mass function of $Y$ is given by
\begin{equation}
p(i)=\binom{m}{i}p^{i}(1-p)^{m-i}, \quad i=0,1,\dots, m. \label{40}
\end{equation}
Here $m$ is a positive integer and $0< p \le 1$.
Then we have
\begin{align}
&E[Y]=mp,\,\,E[e^{Yt}]=\big(1+p(e^{t}-1)\big)^{m}, \label{41}\\
& e_{Y}(t)=\big(1+p(e^{t}-1)\big)^{m}-1, \,\,
\bar{e}_{Y}(t)=\log\Big(1+\frac{1}{p}\big((1+t)^{\frac{1}{m}}-1\big)\Big), \nonumber\\
&\bar{f}_{Y}(t)=m \log \big(1+p(e^{t}-1)\big),\,\, f_{Y}(t)=\log\big(1+\frac{1}{p}(e^{\frac{t}{m}}-1)\big). \nonumber
\end{align}
By using \eqref{41}, we obtain
\begin{align}
\frac{1}{k!}\big(e_{Y}(t)\big)^{k}&=\frac{1}{k!}\Big(\big(1+p(e^{t}-1)
 \big)^{m}-1\Big)^{k} \label{42}\\
&=\frac{1}{k!}\Big(e^{m\log \big(1+p(e^{t}-1)\big)}-1 \Big)^{k} \nonumber\\
&=\sum_{j=k}^{\infty} m^{j} S_{2}(j,k) \frac{1}{j!}\Big(\log \big(1+p(e^{t}-1) \big) \Big)^{j} \nonumber\\
&=\sum_{j=k}^{\infty} m^{j} S_{2}(j,k)\sum_{i=j}^{\infty}S_{1}(i,j)p^{i}\frac{1}{i!}(e^{t}-1)^{i} \nonumber \\
&=\sum_{j=k}^{\infty}\sum_{i=j}^{\infty}m^{j}p^{i}S_{2}(j,k)S_{1}(i,j)\sum_{n=i}^{\infty}S_{2}(n,i)\frac{t^{n}}{n!} \nonumber \\
&=\sum_{n=k}^{\infty}\sum_{j=k}^{n}\sum_{i=j}^{n}m^{j}p^{i}S_{2}(j,k)S_{1}(i,j)S_{2}(n,i)\frac{t^{n}}{n!}. \nonumber
\end{align}
Thus \eqref{42} shows that
\begin{equation}
S_{2}^{Y}(n,k)=\sum_{j=k}^{n}\sum_{i=j}^{n}m^{j}p^{i}S_{2}(j,k)S_{1}(i,j)S_{2}(n,i), \quad (n \ge k). \label{43}
\end{equation}
Again, by utilizing \eqref{41}, we have
\begin{align}
\frac{1}{k!}\big(\bar{e}_{Y}(t)\big)^{k}&=\frac{1}{k!}\bigg(\log \Big(1+\frac{1}{p} \big( (1+t)^{\frac{1}{m}}-1 \big) \Big) \bigg)^{k} \label{44}\\
&=\sum_{j=k}^{\infty}S_{1}(j,k)\frac{1}{p^{j}}\frac{1}{j!}\big((1+t)^{\frac{1}{m}}-1\big)^{j} \nonumber\\
&=\sum_{j=k}^{\infty}S_{1}(j,k)\frac{1}{p^{j}}\frac{1}{j!}\big(e^{\frac{1}{m}\log(1+t)}-1 \big)^{j} \nonumber\\
&=\sum_{j=k}^{\infty}\sum_{i=j}^{\infty}\frac{1}{p^{j}}\frac{1}{m^{i}}S_{1}(j,k)S_{2}(i,j)\frac{1}{i!}\big(\log(1+t) \big)^{i} \nonumber\\
&=\sum_{j=k}^{\infty}\sum_{i=j}^{\infty}\frac{1}{p^{j}}\frac{1}{m^{i}}S_{1}(j,k)S_{2}(i,j)\sum_{n=i}^{\infty}S_{1}(n,i)\frac{t^{n}}{n!} \nonumber \\
&=\sum_{n=k}^{\infty}\sum_{j=k}^{n}\sum_{i=j}^{n}\frac{1}{p^{j}}\frac{1}{m^{i}}S_{1}(j,k)S_{2}(i,j)S_1(n,i) \frac{t^{n}}{n!}. \nonumber
\end{align}
Therefore, from \eqref{44}, we get
\begin{equation}
S_{1}^{Y}(n,k)=\sum_{j=k}^{n}\sum_{i=j}^{n}\frac{1}{p^{j}}\frac{1}{m^{i}}S_{1}(j,k)S_{2}(i,j)S_1(n,i). \label{45}
\end{equation}
We observe from \eqref{32} and \eqref{40} that
\begin{align}
E[e_{\lambda}^{Y}(t)]&=\sum_{i=0}^{m} e_{\lambda}^{i}(t)p(i)=\sum_{i=0}^{m}\binom{m}{i}\big(pe_{\lambda}(t) \big)^{i}(1-p)^{m-i} \label{46} \\
&=\Big(1+p\big(e_{\lambda}(t)-1 \big) \Big)^{m}. \nonumber
\end{align}
From \eqref{46}, we get
\begin{align}
& e_{Y,\lambda}(t)=\big(1+p(e_{\lambda}(t)-1)\big)^{m}-1,\,\,\,
\bar{e}_{Y,\lambda}(t)=\log_{\lambda}\Big(1+\frac{1}{p} \big( (1+t)^{\frac{1}{m}}-1 \big) \Big), \label{47}\\
&\bar{f}_{Y,\lambda}(t)=m\log \big(1+p(e_{\lambda}(t)-1)\big),\,\,\,
f_{Y,\lambda}(t)=\log_{\lambda}\big(1+\frac{1}{p}(e^{\frac{t}{m}}-1) \big). \nonumber
\end{align}
From \eqref{47}, we obtain
\begin{align}
&S_{2,\lambda}^{Y}(n,k)=\sum_{j=k}^{n}\sum_{i=j}^{n}m^{j}p^{i}S_{2}(j,k)S_{1}(i,j)S_{2,\lambda}(n,i), \label{48}\\
&S_{1,\lambda}^{Y}(n,k)=\sum_{j=k}^{n}\sum_{i=j}^{n}\frac{1}{p^{j}}\frac{1}{m^{i}}S_{1,\lambda}(j,k)S_{2}(i,j)S_{1}(n,i). \nonumber
\end{align}
Thus the next theorem follows from \eqref{43}, \eqref{45} and \eqref{48}.
\begin{theorem}
Let $Y$ be the binomial random variable with parameter $(m,p)$, with $p >0$. Then, for\, $n \ge k$, we have
\begin{align*}
&S_{2}^{Y}(n,k)=\sum_{j=k}^{n}\sum_{i=j}^{n}m^{j}p^{i}S_{2}(j,k)S_{1}(i,j)S_{2}(n,i),\\
&S_{1}^{Y}(n,k)=\sum_{j=k}^{n}\sum_{i=j}^{n}\frac{1}{p^{j}}\frac{1}{m^{i}}S_{1}(j,k)S_{2}(i,j)S_1(n,i), \\
&S_{2,\lambda}^{Y}(n,k)=\sum_{j=k}^{n}\sum_{i=j}^{n}m^{j}p^{i}S_{2}(j,k)S_{1}(i,j)S_{2,\lambda}(n,i),\\
&S_{1,\lambda}^{Y}(n,k)=\sum_{j=k}^{n}\sum_{i=j}^{n}\frac{1}{p^{j}}\frac{1}{m^{i}}S_{1,\lambda}(j,k)S_{2}(i,j)S_{1}(n,i).
\end{align*}
\end{theorem}
\vspace{0.1in}
(c) Let $Y$ be the Poisson random variable with parameter $\alpha >0$. Then the probability mass function of $Y$ is given by
\begin{equation}
p(i)=e^{-\alpha}\frac{\alpha^{i}}{i!}, \quad i=0,1,2,\dots. \label{49}
\end{equation}
Then we have
\begin{align}
&E[Y]=\alpha,\,\,E[e^{Yt}]=e^{\alpha(e^{t}-1)},\,\,e_{Y}(t)=e^{\alpha(e^{t}-1)}-1, \label{50}\\
&\bar{e}_{Y}(t)=\log\big(1+\frac{1}{\alpha} \log (1+t) \big),\,\,\bar{f}_{Y}(t)=\alpha(e^{t}-1),\,\, f_{Y}(t)=\log(1+\frac{t}{\alpha}). \nonumber
\end{align}
From \eqref{50}, we get
\begin{equation}
S_{2}^{Y}(n,k)=\sum_{j=k}^{n}\alpha^{j}S_{2}(j,k)S_{2}(n,j),\,\,\, S_{1}^{Y}(n,k)=\sum_{j=k}^{n}\frac{1}{\alpha^{j}}S_{1}(j,k)S_{1}(n,j). \label{51}
\end{equation}
From \eqref{32} and \eqref{49}, we observe that
\begin{equation}
E[e_{\lambda}^{Y}(t)]=\sum_{i=0}^{\infty}e_{\lambda}^{i}(t)p(i)
=e^{-\alpha}\sum_{i=0}^{\infty}\frac{1}{i!} \big(\alpha e_{\lambda}(t) \big)^{i}
=e^{\alpha \big(e_{\lambda}(t)-1 \big)}. \label{52}
\end{equation}
From \eqref{52}, we derive
\begin{align}
&e_{Y,\lambda}(t)=e^{\alpha(e_{\lambda}(t)-1)}-1,\,\,\, \bar{e}_{Y,\lambda}(t)=\log_{\lambda}\big(1+\frac{1}{\alpha}\log(1+t)\big), \label{53}\\
&\bar{f}_{Y,\lambda}(t)=\alpha(e_{\lambda}(t)-1),\,\,\, f_{Y,\lambda}(t)=\log_{\lambda}(1+\frac{t}{\alpha}). \nonumber
\end{align}
From \eqref{53}, we get
\begin{equation}
S_{2,\lambda}^{Y}(n,k)=\sum_{j=k}^{n}\alpha^{j}S_{2}(j,k)S_{2,\lambda}(n,j),\,\,\, S_{1,\lambda}^{Y}(n,k)=\sum_{j=k}^{n}\frac{1}{\alpha^{j}}S_{1,\lambda}(j,k)S_{1}(n,j). \label{54}
\end{equation}
From \eqref{51} and \eqref{54}, the next theorem follows.
\begin{theorem}
Let $Y$ be the Poisson random variable with parameter $\alpha >0$. Then, for\, $n \ge k$, we have
\begin{align*}
&S_{2}^{Y}(n,k)=\sum_{j=k}^{n}\alpha^{j}S_{2}(j,k)S_{2}(n,j),\,\,\, S_{1}^{Y}(n,k)=\sum_{j=k}^{n}\frac{1}{\alpha^{j}}S_{1}(j,k)S_{1}(n,j),\\
&S_{2,\lambda}^{Y}(n,k)=\sum_{j=k}^{n}\alpha^{j}S_{2}(j,k)S_{2,\lambda}(n,j),\,\,\, S_{1,\lambda}^{Y}(n,k)=\sum_{j=k}^{n}\frac{1}{\alpha^{j}}S_{1,\lambda}(j,k)S_{1}(n,j).
\end{align*}
\end{theorem}
\vspace{0.1in}
(d) Let $Y$ be the geometric random variable with parameter $0 <p <1$. Then the probability mass function of $Y$ is given by
\begin{equation}
p(i)=(1-p)^{i-1}p, \quad i=1,2,\dots. \label{55}
\end{equation}
Then we have
\begin{align}
&E[Y]=\frac{1}{p}, \,\, E[e^{Yt}]=\frac{pe^{t}}{1-(1-p)e^{t}}, \label{56}\\
&e_{Y}(t)=\frac{pe^{t}}{1-(1-p)e^{t}}-1,\,\, \bar{e}_{Y}(t)=\log \Big(\frac{1+t}{1+(1-p)t}\Big), \nonumber \\
&\bar{f}_{Y}(t)=\log\Big(\frac{pe^{t}}{1-(1-p)e^{t}} \Big),\,\,
f_{Y}(t)=\log \Big(\frac{e^{t}}{p+(1-p)e^{t}}\Big). \nonumber
\end{align}
Before proceeding further, we observe from \eqref{56} that
\begin{equation}
e_{Y}(t)=\frac{e^{t}-1}{(p-1)e^{t}+1}=\frac{1}{p-1}\bigg(1-\frac{1-\frac{1}{1-p}}{e^{t}-\frac{1}{1-p}} \bigg). \label{57}
\end{equation}
Thus, from \eqref{57}, we have
\begin{align}
\frac{1}{k!}\big(e_{Y}(t) \big)^{k}&=\frac{1}{k!}\frac{1}{(p-1)^{k}}\bigg(1-\frac{1-\frac{1}{1-p}}{e^{t}-\frac{1}{1-p}} \bigg)^{k} \label{58}\\
&=\frac{1}{k!}\frac{1}{(p-1)^{k}}\sum_{j=0}^{k}\binom{k}{j}(-1)^{j}\bigg(\frac{1-\frac{1}{1-p}}{e^{t}-\frac{1}{1-p}} \bigg)^{j} \nonumber \\
&=\frac{1}{k!}\frac{1}{(p-1)^{k}}\sum_{j=0}^{k}\binom{k}{j}(-1)^{j}\sum_{n=0}^{\infty}H_{n}^{(j)}\Big(\frac{1}{1-p}\Big)\frac{t^{n}}{n!} \nonumber \\
&=\sum_{n=0}^{\infty}\frac{1}{k!}\frac{1}{(p-1)^{k}}\sum_{j=0}^{k}\binom{k}{j}(-1)^{j}H_{n}^{(j)}\Big(\frac{1}{1-p}\Big)\frac{t^{n}}{n!}, \nonumber
\end{align}
where, for any nonnegative integer $r$, $H_{n}^{(r)}(u)$ are the Frobenius-Euler numbers of order $r$, given by
\begin{equation}
\bigg(\frac{1-u}{e^{t}-u} \bigg)^{r}=\sum_{n=0}^{\infty}H_{n}^{(r)}(u)\frac{t^{n}}{n!}, \quad(u \ne 1). \label{59}
\end{equation}
Thus we have shown that
\begin{equation}
S_{2}^{Y}(n,k)=\frac{1}{k!}\frac{1}{(p-1)^{k}}\sum_{j=0}^{k}\binom{k}{j}(-1)^{j}H_{n}^{(j)}\Big(\frac{1}{1-p}\Big), \quad \mathrm{for}\,\,\, n \ge k, \label{60}
\end{equation}
and
\begin{equation}
\sum_{j=0}^{k}\binom{k}{j}(-1)^{j}H_{n}^{(j)}\Big(\frac{1}{1-p}\Big)=0, \quad \mathrm{for}\,\,\, 0 \le n < k. \label{61}
\end{equation}
Next, we consider
\begin{align}
\frac{1}{k!}\big(\bar{e}_{Y}(t)\big)^{k}&=\frac{1}{k!}\bigg(\log\Big(\frac{1+t}{1+(1-p)t}\Big)\bigg)^{k} \label{62}\\
&=\frac{1}{k!}\bigg(\log \Big(1+\frac{pt}{1+(1-p)t} \Big) \bigg)^{k} \nonumber \\
&=\sum_{j=k}^{\infty}S_{1}(j,k)\frac{1}{j!}\Big(\frac{pt}{1+(1-p)t} \Big)^{j} \nonumber \\
&=\sum_{j=k}^{\infty}S_{1}(j,k)p^{j}\frac{t^{j}}{j!}\big(1+(1-p)t\big)^{-j} \nonumber \\
&=\sum_{j=k}^{\infty}S_{1}(j,k)p^{j}\frac{t^{j}}{j!}\sum_{m=0}^{\infty}(j+m-1)_{m}(p-1)^{m}\frac{t^{m}}{m!} \nonumber \\
&=\sum_{n=k}^{\infty}\sum_{j=k}^{n}\binom{n}{j}(n-1)_{n-j}p^{j}(p-1)^{n-j}S_{1}(j,k) \frac{t^{n}}{n!}. \nonumber
\end{align}
Thus \eqref{62} shows that
\begin{equation}
S_{1}^{Y}(n,k)=\sum_{j=k}^{n}\binom{n}{j}(n-1)_{n-j}p^{j}(p-1)^{n-j}S_{1}(j,k),\,\,\,\mathrm{for}\,\,\,n \ge k. \label{63}
\end{equation}
We observe from \eqref{32} and \eqref{55} that
\begin{align}
E[e_{\lambda}^{Y}(t)]&=\sum_{i=1}^{\infty}e_{\lambda}^{i}(t)p(i)=pe_{\lambda}(t)\sum_{i=1}^{\infty}\big((1-p)e_{\lambda}(t) \big)^{i-1} \label{64}\\
&=\frac{pe_{\lambda}(t)}{1-(1-p)e_{\lambda}(t)}. \nonumber
\end{align}
From \eqref{64}, we obtain
\begin{align}
&e_{Y,\lambda}(t)=\frac{pe_{\lambda}(t)}{1-(1-p)e_{\lambda}(t)}-1,\,\,\, \bar{e}_{Y,\lambda}(t)=\log_{\lambda} \Big(\frac{1+t}{1+(1-p)t}\Big), \label{65}\\
&\bar{f}_{Y,\lambda}(t)=\log\Big(\frac{pe_{\lambda}(t)}{1-(1-p)e_{\lambda}(t)} \Big),\,\,\,
f_{Y,\lambda}(t)=\log_{\lambda} \Big(\frac{e^{t}}{p+(1-p)e^{t}}\Big). \nonumber
\end{align}
Using \eqref{65} and proceeding just as in \eqref{58} with $e^{t}$ replaced by $e_{\lambda}(t)$, we show that
\begin{equation}
S_{2,\lambda}^{Y}(n,k)=\frac{1}{k!}\frac{1}{(p-1)^{k}}\sum_{j=0}^{k}\binom{k}{j}(-1)^{j}h_{n,\lambda}^{(j)}\Big(\frac{1}{1-p}\Big), \quad \mathrm{for}\,\,\, n \ge k, \label{66}
\end{equation}
and
\begin{equation}
\sum_{j=0}^{k}\binom{k}{j}(-1)^{j}h_{n,\lambda}^{(j)}\Big(\frac{1}{1-p}\Big)=0, \quad \mathrm{for}\,\,\, 0 \le n < k. \label{67}
\end{equation}
Here, for any nonnegative integer $r$, $h_{n,\lambda}^{(r)}(u)$ are the degenerate Frobenius-Euler numbers of order $r$, given by
\begin{equation}
\Big(\frac{1-u}{e_{\lambda}(t)-u} \Big)^{r}=\sum_{n=0}^{\infty} h_{n,\lambda}^{(r)}(u)\frac{t^{n}}{n!}, \quad (u \ne 1). \label{68}
\end{equation}
Again, by using \eqref{65} and proceeding just as in \eqref{62} with $e^{t}$ replaced by $e_{\lambda}(t)$, we get
\begin{equation}
S_{1,\lambda}^{Y}(n,k)=\sum_{j=k}^{n}\binom{n}{j}(n-1)_{n-j}p^{j}(p-1)^{n-j}S_{1,\lambda}(j,k), \quad \mathrm{for}\,\,\, n \ge k. \label{69}
\end{equation}
Now, \eqref{60}, \eqref{63}, \eqref{66} and \eqref{69} altogether give the following theorem.
\begin{theorem}
Let $Y$ be the geometric random variable with parameter $0 <p <1$. Then, for\, $n \ge k$, we have
\begin{align*}
&S_{2}^{Y}(n,k)=\frac{1}{k!}\frac{1}{(p-1)^{k}}\sum_{j=0}^{k}\binom{k}{j}(-1)^{j}H_{n}^{(j)}\Big(\frac{1}{1-p}\Big), \\
&S_{1}^{Y}(n,k)=\sum_{j=k}^{n}\binom{n}{j}(n-1)_{n-j}p^{j}(p-1)^{n-j}S_{1}(j,k),\\
&S_{2,\lambda}^{Y}(n,k)=\frac{1}{k!}\frac{1}{(p-1)^{k}}\sum_{j=0}^{k}\binom{k}{j}(-1)^{j}h_{n,\lambda}^{(j)}\Big(\frac{1}{1-p}\Big), \\
&S_{1,\lambda}^{Y}(n,k)=\sum_{j=k}^{n}\binom{n}{j}(n-1)_{n-j}p^{j}(p-1)^{n-j}S_{1,\lambda}(j,k),
\end{align*}
where $H_{n}^{(j)}(u)$ are the Frobenius-Euler numbers of order $j$ (see \eqref{59}) and $h_{n}^{(j)}(u)$ are the degenerate Frobenius-Euler numbers of order $j$ (see \eqref{68}).
\end{theorem}
\vspace{0.1in}
(e) Let $Y$ be the exponential random variable with parameter $\alpha > 0$. Then the probability density function of $Y$ is given by
\begin{equation}
f(y)=\left\{\begin{array}{ccc}
\alpha e^{-\alpha y}, & \textrm{if \,\,$y \ge 0$,} \\
0, & \textrm{if\,\, $y<0$}. \label{70}
\end{array}\right.
\end{equation}
Then we have
\begin{align}
& E[Y]=\frac{1}{\alpha},\,\,\,E[e^{Yt}]=\frac{\alpha}{\alpha-t}, \,\,\, e_{Y}(t)=\frac{t}{\alpha -t}, \label{71}\\
&\bar{e}_{Y}(t)=\frac{\alpha t}{1+t},\,\,\, \bar{f}_{Y}(t)=\log \big(\frac{\alpha}{\alpha-t} \big),\,\,\, f_{Y}(t)=\alpha (1-e^{-t}). \nonumber
\end{align}
By using \eqref{71}, we have
\begin{align}
\frac{1}{k!}\big(e_{Y}(t)\big)^{k}&=\frac{1}{k!}\Big(\frac{t}{\alpha-t} \Big)^{k}=\frac{1}{k!}\frac{1}{\alpha^{k}}t^{k}\big(1-\frac{t}{\alpha} \big)^{-k} \label{72}\\
&=\frac{1}{k!}\frac{1}{\alpha^{k}}t^{k}\sum_{n=0}^{\infty}(k+n-1)_{n}\frac{1}{n!}\Big(\frac{t}{\alpha} \Big)^{n} \nonumber \\
&=\sum_{n=0}^{\infty}(k+n-1)_{n}\frac{1}{n! k!} \frac{1}{\alpha^{n+k}} t^{n+k} \nonumber \\
&=\sum_{n=k}^{\infty}\binom{n}{k}(n-1)_{n-k}\frac{1}{\alpha^{n}}\frac{t^{n}}{n!}.\nonumber
\end{align}
Thus \eqref{72} shows that
\begin{equation}
S_{2}^{Y}(n,k)=\binom{n}{k}(n-1)_{n-k}\frac{1}{\alpha^{n}}, \quad \mathrm{for}\,\,\, n \ge k. \label{73}
\end{equation}
Next, by using \eqref{71}, we obtain
\begin{align}
\frac{1}{k!}\big(\bar{e}_{Y}(t) \big)^{k}&=\frac{1}{k!}\Big(\frac{\alpha t}{1+t} \Big)
=\frac{\alpha^{k}}{k!}t^{k}(1+t)^{-k} \label{74}\\
&=\frac{\alpha^{k}}{k!}t^{k}\sum_{n=0}^{\infty}(-1)^{n}(k+n-1)_{n} \frac{t^{n}}{n!} \nonumber \\
&=\frac{\alpha^{k}}{k!}t^{k}\sum_{n=k}^{\infty}(-1)^{n-k}(n-1)_{n-k} \frac{t^{n-k}}{(n-k)!} \nonumber \\
&=\sum_{n=k}^{\infty}\binom{n}{k}(-1)^{n-k}(n-1)_{n-k} \alpha^{k} \frac{t^{n}}{n!}. \nonumber
\end{align}
Therefore \eqref{74} yields that
\begin{equation}
S_{1}^{Y}(n,k)=(-1)^{n-k}\binom{n}{k}(n-1)_{n-k}\alpha^{k}, \quad \mathrm{for}\,\,\, n \ge k. \label{75}
\end{equation}
From \eqref{33} and \eqref{70}, we observe
\begin{align}
E[e_{\lambda}^{Y}(t)]&=\int_{-\infty}^{\infty}e_{\lambda}^{y}(t)f(y) dy
=\alpha \int_{0}^{\infty}e^{-y\big(\alpha-\frac{1}{\lambda}\log(1+ \lambda t) \big)} dy \label{76}\\
&=\frac{\alpha}{\alpha-\log e_{\lambda}(t)}, \quad \Big(t <\frac{e^{\lambda \alpha}-1}{\lambda} \Big). \nonumber
\end{align}
From \eqref {76}, we get
\begin{align}
&e_{Y,\lambda}(t)=\frac{\log (1+ \lambda t)}{\alpha \lambda -\log(1+ \lambda t)},\,\,\,
\bar{e}_{Y,\lambda}(t)=\frac{1}{\lambda}\big(e^{\frac{\alpha \lambda t}{1+t}}-1 \big), \label{77}\\
&\bar{f}_{Y,\lambda}(t)=\log \Big(\frac{\alpha \lambda}{\alpha \lambda-\log(1+\lambda t)} \Big),\,\,\, f_{Y,\lambda}(t)=\frac{1}{\lambda} \big( e^{\alpha \lambda (1-e^{-t})} -1 \big). \nonumber
\end{align}
Using \eqref{77} and proceeding just as \eqref{72} with $t$ replaced by $\log(1+\lambda t)$ and $\alpha$ by $\alpha \lambda$, we get
\begin{align}
\frac{1}{k!}\big(e_{Y,\lambda}(t) \big)^{k}&=\frac{1}{k!}\bigg(\frac{\log(1+\lambda t)}{\alpha \lambda-\log(1+\lambda t)} \bigg)^{k} \label{78}\\
&= \sum_{j=k}^{\infty}\binom{j}{k}(j-1)_{j-k}\frac{1}{(\alpha \lambda)^{j}}\frac{1}{j!}\big(\log(1+\lambda t)\big)^{j} \nonumber \\
&= \sum_{j=k}^{\infty}\binom{j}{k}(j-1)_{j-k}\frac{1}{(\alpha \lambda)^{j}}\sum_{n=j}^{\infty}S_{1}(n,j)\lambda^{n}\frac{t^{n}}{n!} \nonumber \\
&=\sum_{n=k}^{\infty}\sum_{j=k}^{n}\binom{j}{k}(j-1)_{j-k}\frac{1}{\alpha^{j}}\lambda^{n-j}S_{1}(n,j)\frac{t^{n}}{n!}. \nonumber
\end{align}
Thus \eqref{78} shows that
\begin{equation}
S_{2,\lambda}^{Y}(n,k)=\sum_{j=k}^{n}\binom{j}{k}(j-1)_{j-k}\frac{1}{\alpha^{j}}\lambda^{n-j}S_{1}(n,j),\quad \mathrm{for}\,\,\, n \ge k. \label{79}
\end{equation}
Using \eqref{77} and \eqref{74} with $\alpha$ replaced by $\alpha \lambda$, we get
\begin{align}
\frac{1}{k!}\big(\bar{e}_{Y,\lambda} \big)^{k}&=\frac{1}{\lambda^{k}}\frac{1}{k!}\Big(e^{\frac{\alpha \lambda t}{1+t}}-1 \Big)^{k}
=\frac{1}{\lambda^{k}}\sum_{j=k}^{\infty}S_{2}(j,k)\frac{1}{j!}\Big(\frac{ \alpha \lambda t}{1+t} \Big)^{j} \label{80}\\
&=\frac{1}{\lambda^{k}}\sum_{j=k}^{\infty}S_{2}(j,k)\sum_{n=j}^{\infty}\binom{n}{j}(-1)^{n-j}(n-1)_{n-j}(\alpha \lambda)^{j}\frac{t^{n}}{n!} \nonumber \\
&=\sum_{n=k}^{\infty}\sum_{j=k}^{n}\binom{n}{j}(-1)^{n-j}(n-1)_{n-j}\alpha^{j} \lambda^{j-k}S_{2}(j,k)\frac{t^{n}}{n!}. \nonumber
\end{align}
Therefore \eqref{80} shows that
\begin{equation}
S_{1,\lambda}^{Y}(n,k)=\sum_{j=k}^{n}\binom{n}{j}(-1)^{n-j}(n-1)_{n-j}\alpha^{j} \lambda^{j-k}S_{2}(j,k), \quad \mathrm{for}\,\,\, n \ge k. \label{81}
\end{equation}
Now, the next result follows \eqref{73}, \eqref{75}, \eqref{79} and \eqref{81}.
\begin{theorem}
 Let $Y$ be the exponential random variable with parameter $\alpha > 0$. Then, for\, $n \ge k$, we have
\begin{align*}
&S_{2}^{Y}(n,k)=\binom{n}{k}(n-1)_{n-k}\frac{1}{\alpha^{n}},\,\,\,S_{1}^{Y}(n,k)=(-1)^{n-k}\binom{n}{k}(n-1)_{n-k}\alpha^{k}, \\
&S_{2,\lambda}^{Y}(n,k)=\sum_{j=k}^{n}\binom{j}{k}(j-1)_{j-k}\frac{1}{\alpha^{j}}\lambda^{n-j}S_{1}(n,j), \\
&S_{1,\lambda}^{Y}(n,k)=\sum_{j=k}^{n}\binom{n}{j}(-1)^{n-j}(n-1)_{n-j}\alpha^{j} \lambda^{j-k}S_{2}(j,k).
\end{align*}
\end{theorem}

\vspace{0.1in}
(f) Let $Y$ be the gamma random variable with parameters $\alpha, \beta >0$. Then the probability density function of $Y$ is given by
\begin{equation}
f(y)=\left\{\begin{array}{ccc}
\beta e^{-\beta y}\frac{(\beta y)^{\alpha-1}}{\Gamma(\alpha)}, & \textrm{if\,\, $y\ge 0$,}\\
0, & \textrm{if \,\,$y<0$}, \label{82}
\end{array}\right.
\end{equation}
Then we have
\begin{align}
&E[Y]=\frac{\alpha}{\beta},\,\,\, E[e^{Yt}]=\Big(\frac{\beta}{\beta-t} \Big)^{\alpha},\,\,\,
e_{Y}(t)=\Big(\frac{\beta}{\beta-t} \Big)^{\alpha}-1, \label{83}\\
&\bar{e}_{Y}(t)=\beta \big(1-(1+t)^{-\frac{1}{\alpha}} \big),\,\,\,\bar{f}_{Y}(t)=\alpha \log \big(\frac{\beta}{\beta-t} \big),\,\,\,f_{Y}(t)=\beta\big(1-e^{-\frac{t}{\alpha}} \big). \nonumber
\end{align}
Using \eqref{83}, we get
\begin{align}
\frac{1}{k!}\big(e_{Y}(t) \big)^{k}&=\frac{1}{k!}\Big(\Big(\frac{\beta}{\beta-t} \Big)^{\alpha}-1 \Big)^{k} \label{84}\\
&=\frac{1}{k!}\sum_{j=0}^{k}\binom{k}{j}(-1)^{k-j}\Big(\frac{1}{1-\frac{t}{\beta}} \Big)^{\alpha j} \nonumber \\
&=\sum_{n=0}^{\infty}\frac{1}{k! \beta^{n}}\sum_{j=0}^{k}\binom{k}{j}(-1)^{k-j}(\alpha j+n-1)_{n} \frac{t^{n}}{n!}. \nonumber
\end{align}
Thus \eqref{84} shows that
\begin{equation}
S_{2}^{Y}(n,k)=\frac{1}{k! \beta^{n}}\sum_{j=0}^{k}\binom{k}{j}(-1)^{k-j}(\alpha j+n-1)_{n},\quad \mathrm{for}\,\,\,n \ge k, \label{85}
\end{equation}
and
\begin{equation}
\sum_{j=0}^{k}\binom{k}{j}(-1)^{k-j}(\alpha j+n-1)_{n}=0,\quad \mathrm{for}\,\,\,0 \le n < k. \label{86}
\end{equation}
Next, by using \eqref{83}, we consider
\begin{align}
\frac{1}{k!}\big(\bar{e}_{Y}(t) \big)^{k}&=\frac{1}{k!}\beta^{k}\big(1-(1+t)^{-\frac{1}{\alpha}} \big)^{k}=\frac{1}{k!}\beta^{k}\sum_{j=0}^{k}\binom{k}{j}(-1)^{j}(1+t)^{-\frac{j}{\alpha}} \label{87} \\
&=\frac{1}{k!}\beta^{k}\sum_{n=0}^{\infty}\big(-\frac{1}{\alpha}\big)^{n}\sum_{j=0}^{k}\binom{k}{j}(-1)^{j}\big(j+(n-1)\alpha)_{n,\alpha}\frac{t^{n}}{n!}. \nonumber
\end{align}
Thus \eqref{87} shows that
\begin{equation}
S_{1}^{Y}(n,k)=\frac{1}{k!}\beta^{k}\big(-\frac{1}{\alpha} \big)^{n}\sum_{j=0}^{k}\binom{k}{j}(-1)^{j}\big(j+(n-1)\alpha \big)_{n,\alpha}, \quad \mathrm{for}\,\,\, n \ge k, \label{88}
\end{equation}
and
\begin{equation}
\sum_{j=0}^{k}\binom{k}{j}(-1)^{j}\big(j+(n-1)\alpha \big)_{n,\alpha}=0, \quad \mathrm{for}\,\,\, 0 \le n < k.  \label{89}
\end{equation}
We observe from \eqref{33} and \eqref{82} that
\begin{align}
E[e_{\lambda}^{Y}(t)]&=\int_{-\infty}^{\infty}e_{\lambda}^{y}f(y) dy
=\frac{1}{\Gamma(\alpha)}\int_{0}^{\infty}e^{-y\big(\beta-\frac{1}{\lambda} \log(1+\lambda t) \big)}(\beta y)^{\alpha} \frac{dy}{y} \label{90}\\
&=\bigg(\frac{\beta}{\beta-\frac{1}{\lambda} \log(1+\lambda t)} \bigg)^{\alpha}, \quad
\Big(t <\frac{e^{\lambda \beta}-1}{\lambda}\Big). \nonumber
\end{align}
From \eqref{90}, we have
\begin{align}
&e_{Y,\lambda}(t)=\bigg(\frac{\beta}{\beta-\frac{1}{\lambda} \log(1+\lambda t)} \bigg)^{\alpha} -1,\,\,\, \bar{e}_{Y,\lambda}(t)=\frac{1}{\lambda} \Big (e^{\lambda \beta \big(1-(1+t)^{-\frac{1}{\alpha}}\big)}-1 \Big), \label{91}\\
&\bar{f}_{Y,\lambda}(t)=\alpha \log\bigg(\frac{\beta}{\beta-\frac{1}{\lambda}\log(1+\lambda t)} \bigg),\,\,\, f_{Y, \lambda}(t)=\frac{1}{\lambda}\Big(e^{\lambda \beta(1-e^{-\frac{t}{\alpha}})}-1 \Big). \nonumber
\end{align}
From \eqref{91} and \eqref{84} with $t$ replaced by $\frac{1}{\lambda}\log(1+\lambda t)$, we obtain
\begin{align}
&\frac{1}{k!}\big(e_{Y,\lambda}(t) \big)^{k}=\frac{1}{k!}\bigg(\Big(\frac{\beta}{\beta-\frac{1}{\lambda} \log(1+\lambda t)}\Big)^{\alpha}-1 \bigg)^{k} \label{92}\\
&=\sum_{l=0}^{\infty}\sum_{j=0}^{k}\frac{1}{k!}\binom{k}{j}(-1)^{k-j}(\alpha j+l-1)_{l}\frac{1}{\beta^{l}}\frac{1}{\lambda^{l}}\frac{1}{l!}\big(\log(1+\lambda t)\big)^{l} \nonumber \\
&=\sum_{l=0}^{\infty}\sum_{j=0}^{k}\frac{1}{k!}\binom{k}{j}(-1)^{k-j}(\alpha j+l-1)_{l}\frac{1}{\beta^{l}}\frac{1}{\lambda^{l}}\sum_{n=l}^{\infty}S_{1}(n,l)\frac{1}{n!}(\lambda t)^{n} \nonumber \\
&=\sum_{n=0}^{\infty}\sum_{l=0}^{n}\sum_{j=0}^{k}(-1)^{k-j}\frac{1}{k!}\binom{k}{j}(\alpha j+l-1)_{l}\frac{1}{\beta^{l}}\lambda^{n-l}S_{1}(n,l) \frac{t^{n}}{n!}. \nonumber
\end{align}
Thus \eqref{92} shows that
\begin{equation}
S_{2,\lambda}^{Y}(n,k)=\sum_{l=0}^{n}\sum_{j=0}^{k}(-1)^{k-j}\frac{1}{k!}\binom{k}{j}(\alpha j+l-1)_{l}\frac{1}{\beta^{l}}\lambda^{n-l}S_{1}(n,l), \quad \mathrm{for}\,\,\, n \ge k, \label{93}
\end{equation}
and
\begin{equation}
\sum_{l=0}^{n}\sum_{j=0}^{k}(-1)^{k-j}\binom{k}{j}(\alpha j+l-1)_{l}\frac{1}{\beta^{l}}\lambda^{n-l}S_{1}(n,l)=0, \quad \mathrm{for}\,\,\, 0 \le n <k. \label{94}
\end{equation}
From \eqref{91} and \eqref{87}, we have
\begin{align}
\frac{1}{k!}\big(\bar{e}_{Y,\lambda}(t) \big)^{k}
&=\frac{1}{\lambda^{k}}\frac{1}{k!} \Big(e^{\lambda \beta (1-(1+t)^{-\frac{1}{\alpha}})} -1 \Big)^{k} \label{95} \\
&=\sum_{l=k}^{\infty}S_{2}(l,k)\lambda^{l-k}\frac{1}{l!}\Big(\beta \big(1-(1+t)^{-\frac{1}{\alpha}} \big) \Big)^{l}  \nonumber \\
&=\sum_{l=k}^{\infty}S_{2}(l,k)\lambda^{l-k}\frac{1}{l!}\beta^{l}\sum_{n=l}^{\infty}\big(-\frac{1}{\alpha} \big)^{n}\sum_{j=0}^{l}\binom{l}{j}(-1)^{j}\big(j+(n-1)\alpha \big)_{n,\alpha}\frac{t^{n}}{n!} \nonumber \\
&=\sum_{n=k}^{\infty}\sum_{l=k}^{n}\sum_{j=0}^{l}S_{2}(l,k)\binom{l}{j}\frac{1}{l!}(-1)^{n-j}\lambda^{l-k}\beta^{l}\frac{1}{\alpha^{n}}\big(j+(n-1)\alpha \big)_{n,\alpha} \frac{t^{n}}{n!}. \nonumber
\end{align}
Therefore \eqref{95} shows that, for\, $n \ge k$,
\begin{equation}
S_{1,\lambda}^{Y}(n,k)=\sum_{l=k}^{n}\sum_{j=0}^{l}S_{2}(l,k)\binom{l}{j}\frac{1}{l!}(-1)^{n-j}\lambda^{l-k}\beta^{l}\frac{1}{\alpha^{n}}\big(j+(n-1)\alpha \big)_{n,\alpha}. \label{96}
\end{equation}
Now, the next result follows from \eqref{85}, \eqref{88}, \eqref{93} and \eqref{96}.
\begin{theorem}
Let $Y$ be the gamma random variable with parameters $\alpha, \beta >0$. Then, for\, $n \ge k$, we have
\begin{align*}
&S_{2}^{Y}(n,k)=\frac{1}{k! \beta^{n}}\sum_{j=0}^{k}\binom{k}{j}(-1)^{k-j}(\alpha j+n-1)_{n}, \\
&S_{1}^{Y}(n,k)=\frac{1}{k!}\beta^{k}\big(-\frac{1}{\alpha} \big)^{n}\sum_{j=0}^{k}\binom{k}{j}(-1)^{j}\big(j+(n-1)\alpha \big)_{n,\alpha}, \\
&S_{2,\lambda}^{Y}(n,k)=\sum_{l=0}^{n}\sum_{j=0}^{k}(-1)^{k-j}\frac{1}{k!}\binom{k}{j}(\alpha j+l-1)_{l}\frac{1}{\beta^{l}}\lambda^{n-l}S_{1}(n,l),\\
&S_{1,\lambda}^{Y}(n,k)=\sum_{l=k}^{n}\sum_{j=0}^{l}S_{2}(l,k)\binom{l}{j}\frac{1}{l!}(-1)^{n-j}\lambda^{l-k}\beta^{l}\frac{1}{\alpha^{n}}\big(j+(n-1)\alpha \big)_{n,\alpha}.
\end{align*}
\end{theorem}
\vspace{0.1in}
(g) Let $Y$ be the normal random variable with parameters $(\mu, \sigma^{2})$, where $\mu \ne 0$. Then the probability density function of $Y$ is given by
\begin{equation}
f(x)=\frac{1}{\sqrt{2 \pi} \sigma}e^{-\frac{(x-\mu)^{2}}{2 \sigma^{2}}},\quad \mathrm{for}\,\,\,-\infty <x < \infty. \label{97}
\end{equation}
This distribution is denoted by $Y \sim N(\mu,\sigma^{2})$. \par
Then we have
\begin{align}
&E[Y]=\mu,\,\,\,E[e^{Yt}]=e^{\frac{\sigma^{2}t^{2}}{2}+\mu t},\,\,\,e_{Y}(t)=e^{\frac{\sigma^{2}t^{2}}{2}+\mu t}-1,\label{98}\\
&\bar{e}_{Y}(t)=-\frac{\mu}{\sigma^{2}}+\frac{\mu}{\sigma^{2}}\Big(1+2\Big(\frac{\sigma}{\mu}\Big)^{2}\log(1+t) \Big)^{\frac{1}{2}}, \nonumber \\
&\bar{f}_{Y}(t)=\frac{\sigma^{2}t^{2}}{2}+\mu t,\,\,\,f_{Y}(t)=-\frac{\mu}{\sigma^{2}}+\Big(\Big(\frac{\mu}{\sigma^{2}}\Big)^{2}+\frac{2t}{\sigma^{2}} \Big)^{\frac{1}{2}}. \nonumber
\end{align}
Here we note from \eqref{98} that
\begin{align*}
\bar{e}_{Y}(t)&=-\frac{\mu}{\sigma^{2}}+\frac{\mu}{\sigma^{2}}\Big(1+2\Big(\frac{\sigma}{\mu}\Big)^{2}\log(1+t) \Big)^{\frac{1}{2}} \\
&=-\frac{\mu}{\sigma^{2}}+\frac{\mu}{\sigma^{2}}\sum_{j=0}^{\infty}\Big(\frac{1}{2}\Big)_{j}2^{j}\Big(\frac{\sigma}{\mu} \Big)^{2j}\frac{1}{j!}\big(\log(1+t)\big)^{j} \\
&=-\frac{\mu}{\sigma^{2}}+\frac{\mu}{\sigma^{2}}\sum_{j=0}^{\infty}\Big(\frac{1}{2}\Big)_{j}2^{j}\Big(\frac{\sigma}{\mu} \Big)^{2j}\sum_{n=j}^{\infty}S_{1}(n,j)\frac{t^{n}}{n!} \\
&=-\frac{\mu}{\sigma^{2}}+\frac{\mu}{\sigma^{2}}\sum_{n=0}^{\infty}\sum_{j=0}^{n}\Big(\frac{1}{2}\Big)_{j}2^{j}\Big(\frac{\sigma}{\mu} \Big)^{2j}S_{1}(n,j)\frac{t^{n}}{n!} \\
&=\frac{1}{\mu}t+\frac{\mu}{\sigma^{2}}\sum_{n=2}^{\infty}\sum_{j=0}^{n}\Big(\frac{1}{2}\Big)_{j}2^{j}\Big(\frac{\sigma}{\mu} \Big)^{2j}S_{1}(n,j)\frac{t^{n}}{n!} .
\end{align*}
From \eqref{98}, we show that
\begin{align}
\frac{1}{k!}\big(e_{Y}(t)\big)^{k}&=\frac{1}{k!}\Big(e^{\frac{\sigma^{2}t^{2}}{2}+\mu t}-1 \Big)^{k} \label{99}\\
&=\sum_{j=k}^{\infty}S_{2}(j,k)\frac{1}{j!}t^{j}\Big(\frac{\sigma^{2}t}{2}+\mu \Big)^{j} \nonumber \\
&=\sum_{j=k}^{\infty}S_{2}(j,k)\frac{1}{j!}t^{j}\sum_{l=0}^{j}\binom{j}{l}\mu^{j-l}\Big(\frac{\sigma^{2}}{2}\Big)^{l}t^{l} \nonumber \\
&=\sum_{j=0}^{\infty}\sum_{l=0}^{\infty}\frac{1}{j!}\binom{j}{l}\mu^{j-l}\Big(\frac{\sigma^{2}}{2}\Big)^{l}S_{2}(j,k)t^{j+l} \nonumber \\
&=\sum_{n=0}^{\infty}\sum_{j=0}^{n}\frac{1}{j!}\binom{j}{n-j}\mu^{2j-n}\Big(\frac{\sigma^{2}}{2} \Big)^{n-j}S_{2}(j,k)t^{n} \nonumber \\
&=\sum_{n=k}^{\infty}\sum_{j=k}^{n}\frac{n!}{j!}\binom{j}{n-j} \mu^{2j-n}  \Big(\frac{\sigma^{2}}{2} \Big)^{n-j} S_{2}(j,k)\frac{t^{n}}{n!}. \nonumber
\end{align}
Thus \eqref{99} shows that
\begin{equation}
S_{2}^{Y}(n,k)=\sum_{j=k}^{n}\frac{n!}{j!}\binom{j}{n-j} \mu^{2j-n}  \Big(\frac{\sigma^{2}}{2} \Big)^{n-j} S_{2}(j,k),\quad \mathrm{for}\,\,\,n \ge k. \label{100}
\end{equation}
From \eqref{98}, we see that
\begin{align}
\frac{1}{k!}\big(\bar{e}_{Y}(t)\big)^{k}&=\frac{1}{k!}\Big(-\frac{\mu}{\sigma^{2}} \Big)^{k}\Big(1-\Big(1+2\Big(\frac{\sigma}{\mu}\Big)^{2}\log(1+t) \Big)^{\frac{1}{2}}\Big)^{k} \label{101}\\
&=\frac{1}{k!}\Big(-\frac{\mu}{\sigma^{2}} \Big)^{k}\sum_{j=0}^{k}\binom{k}{j}(-1)^{j}\Big(1+2\Big(\frac{\sigma}{\mu}\Big)^{2}\log(1+t) \Big)^{\frac{j}{2}} \nonumber
\\
&=\frac{1}{k!}\Big(-\frac{\mu}{\sigma^{2}} \Big)^{k}\sum_{j=0}^{k}\binom{k}{j}(-1)^{j}\sum_{l=0}^{\infty}\Big(\frac{j}{2}\Big)_{l}2^{l}\Big(\frac{\sigma}{\mu}\Big)^{2l}\frac{1}{l!}\big(\log(1+t)\big)^{l} \nonumber\\
&=\frac{1}{k!}\Big(-\frac{\mu}{\sigma^{2}} \Big)^{k}\sum_{j=0}^{k}\binom{k}{j}(-1)^{j}\sum_{l=0}^{\infty}\Big(\frac{j}{2}\Big)_{l}2^{l}\Big(\frac{\sigma}{\mu}\Big)^{2l}\sum_{n=l}^{\infty}S_{1}(n,l)\frac{t^{n}}{n!} \nonumber\\
&=\sum_{n=0}^{\infty}\frac{1}{k!}\Big(-\frac{\mu}{\sigma^{2}}\Big)^{k}\sum_{j=0}^{k}\sum_{l=0}^{n}(-1)^{j}\binom{k}{j}\Big(\frac{j}{2}\Big)_{l}2^{l}\Big(\frac{\sigma}{\mu}\Big)^{2l}S_{1}(n,l) \frac{t^{n}}{n!}. \nonumber
\end{align}
Therefore \eqref{101} shows that
\begin{equation}
S_{1}^{Y}(n,k)=\frac{1}{k!}\Big(-\frac{\mu}{\sigma^{2}}\Big)^{k}\sum_{j=0}^{k}\sum_{l=0}^{n}(-1)^{j}\binom{k}{j}\Big(\frac{j}{2}\Big)_{l}2^{l}\Big(\frac{\sigma}{\mu}\Big)^{2l}S_{1}(n,l), \quad \mathrm{for}\,\,\, n \ge k,\label{102}
\end{equation}
and
\begin{equation}
\sum_{j=0}^{k}\sum_{l=0}^{n}(-1)^{j}\binom{k}{j}\Big(\frac{j}{2}\Big)_{l}2^{l}\Big(\frac{\sigma}{\mu}\Big)^{2l}S_{1}(n,l)=0,\quad \mathrm{for}\,\,\,  0 \le n < k. \label{103}
\end{equation}\par
Let $Z$ be the standard normal distribution, that is $Z \sim N(0,1)$. Then $Y=\sigma Z+\mu \sim N(\mu, \sigma^{2})$ (see \eqref{97}).
Note that
\begin{align}
E[e_{\lambda}^{\sigma Z}(t)]&=\frac{1}{\sqrt{2 \pi}}\int_{-\infty}^{\infty}e_{\lambda}^{\sigma z} e^{-\frac{1}{2}z^{2}} dz \label{104}\\
&=\frac{1}{\sqrt{2 \pi}}\int_{-\infty}^{\infty}e^{-\frac{1}{2}\big(z^{2}-2\frac{\sigma}{\lambda}\log(1+\lambda t)z \big)}dz \nonumber\\
&=\frac{1}{\sqrt{2 \pi}}e^{\frac{1}{2}\big(\sigma\log e_{\lambda}(t)\big)^{2}}\int_{-\infty}^{\infty}e^{-\frac{1}{2}\big(z-\sigma\log e_{\lambda}(t)\big)^{2}} dz \nonumber\\
&=e^{\frac{1}{2}\big(\sigma\log e_{\lambda}(t)\big)^{2}}. \nonumber
\end{align}
Using \eqref{104}, we see that
\begin{equation}
E[e_{\lambda}^{Y}(t)]=e_{\lambda}^{\mu}(t)E[e_{\lambda}^{\sigma Z}(t)]
=e_{\lambda}^{\mu}(t)e^{\frac{1}{2}\big(\sigma\log e_{\lambda}(t)\big)^{2}}
=e^{\mu\log e_{\lambda}(t)+\frac{1}{2}\sigma^{2}\big(\log e_{\lambda}(t)\big)^{2}}\label{105}.
\end{equation}
From \eqref{105}, we get
\begin{align}
&e_{Y,\lambda}(t)=e^{\mu\log e_{\lambda}(t)+\frac{1}{2}\sigma^{2}\big(\log e_{\lambda}(t)\big)^{2}}-1, \label{106}\\
&\bar{e}_{Y,\lambda}(t)=\frac{1}{\lambda}\Big(e^{-\frac{\lambda \mu}{\sigma^{2}}+\frac{\lambda \mu}{\sigma^{2}}\big(1+2\big(\frac{\sigma}{\mu}\big)^{2} \log(1+t) \big)^{\frac{1}{2}}}-1 \Big),\nonumber\\
&\bar{f}_{Y,\lambda}(t)=\mu\log e_{\lambda}(t)+\frac{1}{2}\sigma^{2}\big(\log e_{\lambda}(t)\big)^{2}, \nonumber \\
&f_{Y,\lambda}(t)=\frac{1}{\lambda}\Big(e^{-\frac{\mu \lambda}{\sigma^{2}}+\frac{\mu \lambda}{\sigma^{2}}\big(1+2\big(\frac{\sigma}{\mu}\big)^{2} t \big)^{\frac{1}{2}}}-1 \Big). \nonumber
\end{align}
Here we note that
\begin{align*}
\bar{e}_{Y,\lambda}(t)&=\frac{1}{\lambda}\Big(e^{-\frac{\lambda \mu}{\sigma^{2}}+\frac{\lambda \mu}{\sigma^{2}}\big(1+2\big(\frac{\sigma}{\mu}\big)^{2} \log(1+t) \big)^{\frac{1}{2}}}-1 \Big)\\
&=\frac{1}{\lambda}\sum_{j=1}^{\infty}\frac{1}{j!}\Big(-\frac{\lambda \mu}{\sigma^{2}}\Big)^{j}\Big(1-\Big(1+2\Big(\frac{\sigma}{\mu}\Big)^{2} \log (1+t) \Big)^{\frac{1}{2}} \Big)^{j} \\
&=\frac{1}{\lambda}\sum_{j=1}^{\infty}\frac{1}{j!}\Big(-\frac{\lambda \mu}{\sigma^{2}}\Big)^{j}\sum_{l=0}^{j}\binom{j}{l}(-1)^{l}\Big(1+2\Big(\frac{\sigma}{\mu}\Big)^{2} \log (1+t) \Big)^{\frac{l}{2}} \\
&=\frac{1}{\lambda}\sum_{j=1}^{\infty}\frac{1}{j!}\Big(-\frac{\lambda \mu}{\sigma^{2}}\Big)^{j}\sum_{l=0}^{j}\binom{j}{l}(-1)^{l}\sum_{m=0}^{\infty}\Big(\frac{l}{2}\Big)_{m}2^{m}\Big(\frac{\sigma}{\mu}\Big)^{2m}\frac{1}{m!}\big(\log(1+t)\big)^{m}\\
&=\frac{1}{\lambda}\sum_{j=1}^{\infty}\frac{1}{j!}\Big(-\frac{\lambda \mu}{\sigma^{2}}\Big)^{j}\sum_{l=1}^{j}\binom{j}{l}(-1)^{l}\sum_{m=0}^{\infty}\Big(\frac{l}{2}\Big)_{m}2^{m}\Big(\frac{\sigma}{\mu}\Big)^{2m}\\
& \qquad\qquad\qquad \times \sum_{n=m}^{\infty}S_{1}(n,m)\frac{t^{n}}{n!}+\frac{1}{\lambda}\big(e^{-\frac{\lambda \mu}{\sigma^{2}}}-1 \big)\\
&=\frac{1}{\lambda}\sum_{n=1}^{\infty}\sum_{m=0}^{n}\sum_{j=1}^{\infty}\sum_{l=1}^{j}(-1)^{j+l}\frac{1}{j!}\binom{j}{l}\Big(\frac{\lambda \mu}{\sigma^{2}}\Big)^{j}\Big(\frac{l}{2}\Big)_{m}2^{m}\Big(\frac{\sigma}{\mu}\Big)^{2m}S_{1}(n,m)\frac{t^{n}}{n!}.
\end{align*}
From \eqref{106}, we have
\begin{align}
\frac{1}{k!}\big(e_{Y,\lambda}(t)\big)^{k}&=\frac{1}{k!}\Big(e^{\mu\log e_{\lambda}(t)+\frac{1}{2}\sigma^{2}\big(\log e_{\lambda}(t)\big)^{2}}-1\Big)^{k}\label{107}\\
&=\sum_{j=k}^{\infty}S_{2}(j,k)\frac{1}{j!}\big(\log e_{\lambda}(t) \big)^{j}\sum_{l=0}^{j}\binom{j}{l}\mu^{j-l}\Big(\frac{\sigma^{2}}{2} \Big)^{l}\big(\log e_{\lambda}(t) \big)^{l} \nonumber \\
&=\sum_{j=k}^{\infty}\sum_{l=0}^{j}\frac{1}{j!}\binom{j}{l}\mu^{j-l}\Big(\frac{\sigma^{2}}{2} \Big)^{l} S_{2}(j,k)\Big(\frac{\log(1+\lambda t)}{\lambda} \Big)^{j+l} \nonumber\\
&=\sum_{j=0}^{\infty}\sum_{l=0}^{\infty}\frac{1}{j!}\binom{j}{l}\mu^{j-l}\Big(\frac{\sigma^{2}}{2} \Big)^{l} S_{2}(j,k)\Big(\frac{\log(1+\lambda t)}{\lambda} \Big)^{j+l} \nonumber\\
&=\sum_{m=0}^{\infty}\sum_{j=0}^{m}\frac{1}{j!}\binom{j}{m-j}\mu^{2j-m}\Big(\frac{\sigma^{2}}{2}\Big)^{m-j}S_{2}(j,k)\lambda^{-m}m! \frac{1}{m!}\big(\log(1+\lambda t)\big)^{m}\nonumber\\
&=\sum_{m=0}^{\infty}\sum_{j=0}^{m}\frac{1}{j!}\binom{j}{m-j}\mu^{2j-m}\Big(\frac{\sigma^{2}}{2}\Big)^{m-j}S_{2}(j,k)\lambda^{-m}m!\sum_{n=m}^{\infty}S_{1}(n,m)\frac{1}{n!}\lambda^{n}t^{n} \nonumber\\
&=\sum_{n=k}^{\infty}\sum_{m=k}^{n}\sum_{j=k}^{m}\frac{m!}{j!}\binom{j}{m-j}\mu^{2j-m}\Big(\frac{\sigma^{2}}{2}\Big)^{m-j}\lambda^{n-m}S_{2}(j,k)S_{1}(n,m)\frac{t^{n}}{n!}. \nonumber
\end{align}
Thus \eqref{107} shows that, for $n \ge k$,
\begin{equation}
S_{2,\lambda}^{Y}(n,k)=\sum_{m=k}^{n}\sum_{j=k}^{m}\frac{m!}{j!}\binom{j}{m-j}\mu^{2j-m}\Big(\frac{\sigma^{2}}{2}\Big)^{m-j}\lambda^{n-m}S_{2}(j,k)S_{1}(n,m). \label{108}
\end{equation}
Next, we consider
\begin{align}
&\frac{1}{k!}\big(\bar{e}_{Y,\lambda}(t)\big)^{k}=\frac{1}{\lambda^{k}}\frac{1}{k!}\Big(e^{-\frac{\lambda \mu}{\sigma^{2}}+\frac{\lambda \mu}{\sigma^{2}}\big(1+2\big(\frac{\sigma}{\mu}\big)^{2} \log(1+t) \big)^{\frac{1}{2}}}-1 \Big)^{k} \label{109}\\
&=\frac{1}{\lambda^{k}}\sum_{j=k}^{\infty}S_{2}(j,k)\frac{1}{j!}\Big(-\frac{\lambda \mu}{\sigma^{2}} \Big)^{j}\Big(1-\Big(1+2\Big(\frac{\sigma}{\mu}\Big)^{2} \log(1+t) \Big)^{\frac{1}{2}} \Big)^{j} \nonumber\\
&=\frac{1}{\lambda^{k}}\sum_{j=k}^{\infty}S_{2}(j,k)\frac{1}{j!}\Big(-\frac{\lambda \mu}{\sigma^{2}} \Big)^{j}\sum_{l=0}^{j}\binom{j}{l}(-1)^{l}\Big(1+2\Big(\frac{\sigma}{\mu}\Big)^{2} \log(1+t) \Big)^{\frac{l}{2}} \nonumber\\
&=\frac{1}{\lambda^{k}}\sum_{j=k}^{\infty}S_{2}(j,k)\frac{1}{j!}\Big(-\frac{\lambda \mu}{\sigma^{2}} \Big)^{j}\sum_{l=0}^{j}\binom{j}{l}(-1)^{l}\sum_{m=0}^{\infty}\Big(\frac{l}{2}\Big)_{m}2^{m}\Big(\frac{\sigma}{\mu} \Big)^{2m}\frac{1}{m!}\big(\log(1+t)\big)^{m} \nonumber\\
&=\frac{1}{\lambda^{k}}\sum_{j=k}^{\infty}S_{2}(j,k)\frac{1}{j!}\Big(-\frac{\lambda \mu}{\sigma^{2}} \Big)^{j}\sum_{l=0}^{j}\binom{j}{l}(-1)^{l}\sum_{m=0}^{\infty}\Big(\frac{l}{2}\Big)_{m}2^{m}\Big(\frac{\sigma}{\mu} \Big)^{2m}\sum_{n=m}^{\infty}S_{1}(n,m)\frac{t^{n}}{n!} \nonumber\\
&=\frac{1}{\lambda^{k}}\sum_{n=0}^{\infty}\sum_{m=0}^{n}\sum_{j=k}^{\infty}\sum_{l=0}^{j}(-1)^{j+l}\frac{1}{j!}\Big(\frac{\lambda \mu}{\sigma^{2}} \Big )^{j}\Big(\frac{l}{2}\Big)_{m}2^{m}\Big(\frac{\sigma}{\mu} \Big)^{2m}\binom{j}{l}S_{2}(j,k) S_{1}(n,m) \frac{t^{n}}{n!}. \nonumber
\end{align}
Therefore \eqref{109} shows that, for \,$n \ge k$,
\begin{equation}
S_{1,\lambda}^{Y}(n,k)=\frac{1}{\lambda^{k}}\sum_{m=0}^{n}\sum_{j=k}^{\infty}\sum_{l=0}^{j}(-1)^{j+l}\frac{1}{j!}\Big(\frac{\lambda \mu}{\sigma^{2}} \Big )^{j}\Big(\frac{l}{2}\Big)_{m}2^{m}\Big(\frac{\sigma}{\mu} \Big)^{2m}\binom{j}{l}S_{2}(j,k) S_{1}(n,m), \label{110}
\end{equation}
and, for \, $0 \le n < k$.
\begin{equation}
\sum_{m=0}^{n}\sum_{j=k}^{\infty}\sum_{l=0}^{j}(-1)^{j+l}\frac{1}{j!}\Big(\frac{\lambda \mu}{\sigma^{2}} \Big )^{j}\Big(\frac{l}{2}\Big)_{m}2^{m}\Big(\frac{\sigma}{\mu} \Big)^{2m}\binom{j}{l}S_{2}(j,k) S_{1}(n,m)=0. \label{111}
\end{equation}
From \eqref{100}, \eqref{102}, \eqref{108} and \eqref{110}, we get the following theorem.
\begin{theorem}
Let $Y$ be the normal random variable with parameters $(\mu, \sigma^{2})$, where $\mu \ne 0$. Then, for $n \ge k$, we have
\begin{align*}
&S_{2}^{Y}(n,k)=\sum_{j=k}^{n}\frac{n!}{j!}\binom{j}{n-j} \mu^{2j-n}  \Big(\frac{\sigma^{2}}{2} \Big)^{n-j} S_{2}(j,k),\\
&S_{1}^{Y}(n,k)=\frac{1}{k!}\Big(-\frac{\mu}{\sigma^{2}}\Big)^{k}\sum_{j=0}^{k}\sum_{l=0}^{n}(-1)^{j}\binom{k}{j}\Big(\frac{j}{2}\Big)_{l}2^{l}\Big(\frac{\sigma}{\mu}\Big)^{2l}S_{1}(n,l), \\
&S_{2,\lambda}^{Y}(n,k)=\sum_{m=k}^{n}\sum_{j=k}^{m}\frac{m!}{j!}\binom{j}{m-j}\mu^{2j-m}\Big(\frac{\sigma^{2}}{2}\Big)^{m-j}\lambda^{n-m}S_{2}(j,k)S_{1}(n,m)
,\\
&S_{1,\lambda}^{Y}(n,k)=\frac{1}{\lambda^{k}}\sum_{m=0}^{n}\sum_{j=k}^{\infty}\sum_{l=0}^{j}(-1)^{j+l}\frac{1}{j!}\Big(\frac{\lambda \mu}{\sigma^{2}} \Big )^{j}\Big(\frac{l}{2}\Big)_{m}2^{m}\Big(\frac{\sigma}{\mu} \Big)^{2m}\binom{j}{l}S_{2}(j,k) S_{1}(n,m).
\end{align*}
\end{theorem}
\vspace{0.1in}
(h) Let $Y$ be the uniform random variable over $(a,b)$, with $a+b \ne 0$. Then the probability density function of $Y$ is given by
\begin{equation}
f(y)=\left\{\begin{array}{ccc}
\frac{1}{b-a}, & \textrm{for\,\, $a < y < b$,}\\
0, & \textrm{otherwise}. \label{112}
\end{array}\right.
\end{equation}
Then we have
\begin{equation}
E[Y]=\frac{1}{2}(a+b),\,\,\, E[e^{Yt}]=\frac{e^{bt}-e^{at}}{(b-a)t},\,\,\,e_{Y}(t)=\frac{e^{bt}-e^{at}}{(b-a)t}-1,\,\,\,\bar{f}_{Y}(t)=\log\Big(\frac{e^{bt}-e^{at}}{(b-a)t}\Big). \label{113}
\end{equation}
From \eqref{113}, we have
\begin{align}
\frac{1}{k!}\big(e_{Y}(t)\big)^{k}&=\frac{1}{k!}\Big(\frac{e^{at}(e^{(b-a)t}-1)}{(b-a)t}-1 \Big)^{k} \label{114}\\
&=\frac{1}{k!}\sum_{j=0}^{k}\binom{k}{j}(-1)^{k-j}e^{ajt}\big((b-a)t\big)^{-j}j!\frac{1}{j!}\big(e^{(b-a)t}-1\big)^{j} \nonumber\\
&=\frac{1}{k!}\sum_{j=0}^{k}\binom{k}{j}(-1)^{k-j}e^{ajt}j!\sum_{l=j}^{\infty}S_{2}(l,j)\frac{1}{l!}\big((b-a)t\big)^{l-j} \nonumber\\
&=\frac{1}{k!}\sum_{j=0}^{k}\binom{k}{j}(-1)^{k-j}e^{ajt}\sum_{l=0}^{\infty}S_{2}(l+j,j)\frac{j!}{(l+j)!}\big((b-a)t\big)^{l} \nonumber\\
&=\frac{1}{k!}\sum_{j=0}^{k}\binom{k}{j}(-1)^{k-j}\sum_{m=0}^{\infty}(aj)^{m}\frac{t^{m}}{m!}\sum_{l=0}^{\infty}S_{2}(l+j,j)\binom{l+j}{j}^{-1}(b-a)^{l}\frac{t^{l}}{l!} \nonumber\\
&=\frac{1}{k!}\sum_{j=0}^{k}\binom{k}{j}(-1)^{k-j}\sum_{n=0}^{\infty}\sum_{l=0}^{n}\binom{n}{l}\binom{l+j}{j}^{-1}(aj)^{n-l}(b-a)^{l}S_{2}(l+j,j)\frac{t^{n}}{n!} \nonumber\\
&=\sum_{n=0}^{\infty}\frac{1}{k!}\sum_{j=0}^{k}\sum_{l=0}^{n}\binom{n}{l}\binom{k}{j}\binom{l+j}{j}^{-1}(-1)^{k-j}(aj)^{n-l}(b-a)^{l}S_{2}(l+j,j)\frac{t^{n}}{n!}. \nonumber
\end{align}
Thus \eqref{114} shows that, for \, $n \ge k$,
\begin{equation}
S_{2}^{Y}(n,k)=\frac{1}{k!}\sum_{j=0}^{k}\sum_{l=0}^{n}\binom{n}{l}\binom{k}{j}\binom{l+j}{j}^{-1}(-1)^{k-j}(aj)^{n-l}(b-a)^{l}S_{2}(l+j,j), \label{115}
\end{equation}
and, for \, $0 \le n <k$,
\begin{equation}
\sum_{j=0}^{k}\sum_{l=0}^{n}\binom{n}{l}\binom{k}{j}\binom{l+j}{j}^{-1}(-1)^{k-j}(aj)^{n-l}(b-a)^{l}S_{2}(l+j,j)=0. \label{116}
\end{equation}
In the special case of $a=0$, we have
\begin{equation}
S_{2}^{Y}(n,k)=\frac{1}{k!}b^n\sum_{j=0}^{k}\binom{k}{j}\binom{n+j}{j}^{-1}(-1)^{k-j}S_{2}(n+j,j),\quad \mathrm{for}\,\,\, n \ge k,  \label{117}
\end{equation}
and
\begin{equation}
\sum_{j=0}^{k}\binom{k}{j}\binom{n+j}{j}^{-1}(-1)^{k-j}S_{2}(n+j,j)=0, \quad \mathrm{for}\,\,\, 0 \le n <k. \label{118}
\end{equation}
Next, we consider the degenerate case. First, we observe that
\begin{align}
E[e_{\lambda}^{Y}(t)]&=\int_{-\infty}^{\infty}e_{\lambda}^{y}(t)f(y) dy=\frac{1}{b-a}\int_{a}^{b}e_{\lambda}^{y}(t) dy \label{119}\\
&=\frac{1}{b-a}\int_{a}^{b}e^{y\log e_{\lambda}(t)} dy =\frac{e_{\lambda}^{b}(t)-e_{\lambda}^{a}(t)}{(b-a)\log e_{\lambda}(t)}.\nonumber
\end{align}
From \eqref{119}, we get
\begin{equation}
e_{Y,\lambda}(t)=\frac{e_{\lambda}^{b}(t)-e_{\lambda}^{a}(t)}{(b-a)\log e_{\lambda}(t)}-1,\,\,\,\bar{f}_{Y,\lambda}(t)=\log \bigg(\frac{e_{\lambda}^{b}(t)-e_{\lambda}^{a}(t)}{(b-a)\log e_{\lambda}(t)}\bigg). \label{120}
\end{equation}
From \eqref{120} and using \eqref{114} with $t$ replaced by $\log e_{\lambda}(t)=\frac{1}{\lambda}\log(1+\lambda t)$, we have
\begin{align}
\frac{1}{k!}\big(e_{Y,\lambda}(t) \big)^{k}&=\frac{1}{k!}\Bigg(\frac{e_{\lambda}^{b}(t)-e_{\lambda}^{a}(t)}{(b-a)\log e_{\lambda}(t)}-1 \bigg)^{k} \label{121}\\
&=\sum_{m=0}^{\infty}\sum_{j=0}^{k}\sum_{l=0}^{m}\frac{1}{k!}\binom{m}{l}\binom{k}{j}\binom{l+j}{j}^{-1}(-1)^{k-j}(aj)^{m-l}(b-a)^{l} \nonumber\\
&\quad\quad\quad \times S_{2}(l+j,j)\frac{1}{\lambda^{m}}\frac{1}{m!}\big(\log(1+\lambda t \big)^{m} \nonumber\\
&=\sum_{m=0}^{\infty}\sum_{j=0}^{k}\sum_{l=0}^{m}\frac{1}{k!}\binom{m}{l}\binom{k}{j}\binom{l+j}{j}^{-1}(-1)^{k-j}(aj)^{m-l}(b-a)^{l} \nonumber \\
&\quad\quad\quad \times S_{2}(l+j,j)\frac{1}{\lambda^{m}}\sum_{n=m}^{\infty}S_{1}(n,m)\frac{\lambda^{n}t^{n}}{n!} \nonumber \\
&=\sum_{n=0}^{\infty}\sum_{m=0}^{n}\sum_{j=0}^{k}\sum_{l=0}^{m}\frac{1}{k!}\binom{m}{l}\binom{k}{j}\binom{l+j}{j}^{-1}(-1)^{k-j}(aj)^{m-l}(b-a)^{l} \nonumber \\
&\quad\quad\quad \times \lambda^{n-m}S_{2}(l+j,j)S_{1}(n,m)\frac{t^{n}}{n!}. \nonumber
\end{align}
Thus \eqref{121} shows that, for $n \ge k$,
\begin{align}
S_{2,\lambda}^{Y}(n,k)&=\frac{1}{k!}\sum_{m=0}^{n}\sum_{j=0}^{k}\sum_{l=0}^{m}\binom{m}{l}\binom{k}{j}\binom{l+j}{j}^{-1}(-1)^{k-j}(aj)^{m-l}(b-a)^{l} \label{122}\\
&\quad\quad\quad \times \lambda^{n-m}S_{2}(l+j,j)S_{1}(n,m), \nonumber
\end{align}
and, for  $0 \le n < k$,
\begin{align}
&\sum_{m=0}^{n}\sum_{j=0}^{k}\sum_{l=0}^{m}\binom{m}{l}\binom{k}{j}\binom{l+j}{j}^{-1}(-1)^{k-j}(aj)^{m-l}(b-a)^{l} \label{123}\\
&\quad\quad\quad \times \lambda^{n-m}S_{2}(l+j,j)S_{1}(n,m)=0. \nonumber
\end{align}
In the special case of $a=0$, we get,\, for $ n \ge k$,
\begin{equation}
S_{2,\lambda}^{Y}(n,k)=\frac{1}{k!}\sum_{m=0}^{n}\sum_{j=0}^{k}\binom{k}{j}\binom{m+j}{j}^{-1}(-1)^{k-j}b^{m}\lambda^{n-m}S_{2}(m+j,j)S_{1}(n,m), \label{124}
\end{equation}
and, \, for $0 \le n <k$.
\begin{equation}
\sum_{m=0}^{n}\sum_{j=0}^{k}\binom{k}{j}\binom{m+j}{j}^{-1}(-1)^{k-j}b^{m}\lambda^{n-m}S_{2}(m+j,j)S_{1}(n,m)=0. \label{125}
\end{equation}
From \eqref{115} and \eqref{122}, we obtain the following theorem.
\begin{theorem}
Let $Y$ be the uniform random variable over $(a,b)$, with $a+b \ne 0$. For $n \ge k$, we have
\begin{align*}
&S_{2}^{Y}(n,k)=\frac{1}{k!}\sum_{j=0}^{k}\sum_{l=0}^{n}\binom{n}{l}\binom{k}{j}\binom{l+j}{j}^{-1}(-1)^{k-j}(aj)^{n-l}(b-a)^{l}S_{2}(l+j,j), \\
&S_{2,\lambda}^{Y}(n,k)=\frac{1}{k!}\sum_{m=0}^{n}\sum_{j=0}^{k}\sum_{l=0}^{m}\binom{m}{l}\binom{k}{j}\binom{l+j}{j}^{-1}(-1)^{k-j}(aj)^{m-l}(b-a)^{l} \\
&\qquad\qquad\qquad \times \lambda^{n-m}S_{2}(l+j,j)S_{1}(n,m).
\end{align*}
\end{theorem}
In the special case of $a=0$, from \eqref{117} and \eqref{124}, we have
\begin{corollary}
Let $Y$ be the uniform random variable over $(0,b)$, with $b >0$. For $ n \ge k$, we have
\begin{align*}
&S_{2}^{Y}(n,k)=\frac{1}{k!}b^n\sum_{j=0}^{k}\binom{k}{j}\binom{n+j}{j}^{-1}(-1)^{k-j}S_{2}(n+j,j), \\
&S_{2,\lambda}^{Y}(n,k)=\frac{1}{k!}\sum_{m=0}^{n}\sum_{j=0}^{k}\binom{k}{j}\binom{m+j}{j}^{-1}(-1)^{k-j}b^{m}\lambda^{n-m}S_{2}(m+j,j)S_{1}(n,m).
\end{align*}
\end{corollary}

Unlike the cases of $(a)-(g)$, we were not able to find some expressions for $\bar{e}_{Y}(t)$ and $\bar{e}_{Y,\lambda}(t)$ in the case of (h). Consequently, we were not able to determine $S_{1}^{Y}(n,k),\,\,\,f_{Y}(t),\,\,\,S_{1,\lambda}^{Y}(n,k)$, and $f_{Y,\lambda}(t)$. We leave this as an open question to the reader.

\section{\bf Some applications}
By using Proposition 1.1 (a), Proposition 1.2 (a), and Theorem 2.4, we get the following proposition.
\begin{proposition}
Let $p$ be a real number satisfying $0 < p <1$. Then we have the following orthogonality relations:
\begin{align*}
&\sum_{k=l}^{n}\sum_{j=0}^{k}\sum_{m=l}^{k}(-1)^{j}\frac{1}{k!}\binom{k}{j}\binom{k}{m}(k-1)_{k-m}\Big(\frac{p}{p-1}\Big)^{m}H_{n}^{(j)}\big(\frac{1}{1-p} \big)S_{1}(m,l)=\delta_{n,l}, \\
&\sum_{k=l}^{n}\sum_{j=k}^{n}\sum_{m=0}^{l}(-1)^{m}\frac{1}{l!}\binom{n}{j}\binom{l}{m}(n-1)_{n-j}\Big(\frac{p}{p-1}\Big)^{j}(p-1)^{n-l}\\
&\qquad \qquad  \times S_{1}(j,k)H_{k}^{(m)}\big(\frac{1}{1-p}\big)=\delta_{n,l},\\
&\sum_{k=l}^{n}\sum_{j=0}^{k}\sum_{m=l}^{k}(-1)^{j}\frac{1}{k!}\binom{k}{j}\binom{k}{m}(k-1)_{k-m}\Big(\frac{p}{p-1}\Big)^{m}h_{n,\lambda}^{(j)}\big(\frac{1}{1-p} \big)S_{1,\lambda}(m,l)=\delta_{n,l}, \\
&\sum_{k=l}^{n}\sum_{j=k}^{n}\sum_{m=0}^{l}(-1)^{m}\frac{1}{l!}\binom{n}{j}\binom{l}{m}(n-1)_{n-j}\Big(\frac{p}{p-1}\Big)^{j}(p-1)^{n-l}\\
&\qquad \qquad  \times S_{1,\lambda}(j,k)h_{k,\lambda}^{(m)}\big(\frac{1}{1-p}\big)=\delta_{n,l},
\end{align*}
where $H_{n}^{(j)}(u)$ are the Frobenius-Euler numbers of order $j$ (see \eqref{59}) and $h_{n}^{(j)}(u)$ are the degenerate Frobenius-Euler numbers of order $j$ (see \eqref{68}).
\end{proposition}
\begin{remark}
The interested reader can obtain other orthogonality relations by using the results in Theorems 2.1-2.3 and Theorems 2.5-2.7. We note here that many important orthogonality and inverse relations can be found in Chapters 2 and 3 of [25] and Chapter 5 of [26].
\end{remark}
From the definition of the probabilistic Euler polynomials $\mathcal{E}_{n}^{Y}(x)$ in \eqref{26}, we see that $\mathcal{E}_{n}^{Y}(x)$ is Sheffer for the pair $\big(\frac{e^{t}+1}{2}, f_{Y}(t) \big)$, with $\bar{f}_{Y}(t)=\log E[e^{Yt}]$. For the basics about umbral calculus, one refers to [26]. However, no knowledge about umbral calculus is needed in this paper. We observe that
\begin{align}
\sum_{n=0}^{\infty}\big(\mathcal{E}_{n}^{Y}(x+1)+\mathcal{E}_{n}^{Y}(x)\big)\frac{t^{n}}{n!}&=2\big(1+E[e^{Yt}]-1 \big)^{x} \label{126}\\
&=2\sum_{k=0}^{\infty}(x)_{k}\frac{1}{k!}\big(E[e^{Yt}]-1 \big)^{k} \nonumber \\
&=2 \sum_{n=0}^{\infty}\sum_{k=0}^{n}(x)_{k}S_{2}^{Y}(n,k)\frac{t^{n}}{n!}. \nonumber
\end{align}
Thus \eqref{126} shows that
\begin{equation}
\big(\mathcal{E}_{n}^{Y}(x+1)+\mathcal{E}_{n}^{Y}(x)\big)=2\sum_{k=0}^{n}S_{2}^{Y}(n,k)(x)_{k}. \label{127}
\end{equation} \par
Let $q(x)$ be any polynomial of degree $n$ with complex coefficients, and let
\begin{equation}
q(x)=\sum_{k=0}^{n}a_{k}\mathcal{E}_{k}^{Y}(x). \label{128}
\end{equation}
Now, we consider
\begin{equation}
a(x)=q(x+1)+q(x).\label{129}
\end{equation}
Then, by using \eqref{127}-\eqref{129}, we have
\begin{align}
a(x)&=\sum_{k=0}^{n}a_{k}\big(\mathcal{E}_{k}^{Y}(x+1)+\mathcal{E}_{k}^{Y}(x)\big) \label{130}\\
&=2\sum_{k=0}^{n}\sum_{j=0}^{k}a_{k}S_{2}^{Y}(k,j)(x)_{j}. \nonumber
\end{align}
From \eqref{130} and by using \eqref{29}-\eqref{31}, we obtain
\begin{equation}
\Delta^{r}a(x)=2\sum_{k=r}^{n}\sum_{j=r}^{k}a_{k}S_{2}^{Y}(k,j)(j)_{r}(x)_{j-r}. \label{131}
\end{equation}
Hence, from \eqref{131}, we get
\begin{equation}
\frac{1}{2r!} \Delta^{r}a(x)|_{x=0}=\sum_{k=r}^{n}S_{2}^{Y}(k,r)a_{k}. \label{132}
\end{equation}
From \eqref{129} and \eqref{132}, and by invoking the inversion relation in Proposition 1.1, we have
\begin{align}
a_{r}&=\frac{1}{2}\sum_{j=r}^{n}S_{1}^{Y}(j,r)\frac{1}{j!}\Delta^{j}a(x)|_{x=0} \label{133}\\
&=\frac{1}{2}\sum_{j=r}^{n}S_{1}^{Y}(j,r)\frac{1}{j!}\big(\Delta^{j}q(1)+\Delta^{j}q(0) \big). \label{134}
\end{align}
From \eqref{133} and \eqref{28}, we obtain the alternative expressions given by
\begin{align}
a_{r}&=\frac{1}{2}\sum_{j=r}^{n}\sum_{i=0}^{j}(-1)^{j-i}\frac{1}{j!}\binom{j}{i}S_{1}^{Y}(j,r)\big(q(i+1)+q(i)\big) \label{135}\\
&=\frac{1}{2}\sum_{j=r}^{n}\sum_{l=j}^{n}S_{1}^{Y}(j,r)S_{2}(l,j)\frac{1}{l!}\big(q^{(l)}(1)+q^{(l)}(0) \big). \nonumber
\end{align}
Thus, from \eqref{134} and \eqref{135}, we obtain the following theorem.
\begin{theorem}
Let $q(x)$ be any polynomial of degree $n$ with complex coefficients, and let
\begin{equation*}
q(x)=\sum_{r=0}^{n}a_{r}\mathcal{E}_{r}^{Y}(x).
\end{equation*}
Then the coefficients are determined by
\begin{align*}
a_{r}&=\frac{1}{2}\sum_{j=r}^{n}S_{1}^{Y}(j,r)\frac{1}{j!}\big(\Delta^{j}q(1)+\Delta^{j}q(0) \big) \\
&=\frac{1}{2}\sum_{j=r}^{n}\sum_{i=0}^{j}(-1)^{j-i}\frac{1}{j!}\binom{j}{i}S_{1}^{Y}(j,r)\big(q(i+1)+q(i)\big) \\
&=\frac{1}{2}\sum_{j=r}^{n}\sum_{l=j}^{n}S_{1}^{Y}(j,r)S_{2}(l,j)\frac{1}{l!}\big(q^{(l)}(1)+q^{(l)}(0) \big).
\end{align*}
\end{theorem}
Proceeding analogously to the above, we can show the next theorem.
\begin{theorem}
Let $q(x)$ be any polynomial of degree $n$ with complex coefficients, and let
\begin{equation*}
q(x)=\sum_{r=0}^{n}a_{r}\mathcal{E}_{r,\lambda}^{Y}(x).
\end{equation*}
Then the coefficients are determined by
\begin{align*}
a_{r}&=\frac{1}{2}\sum_{j=r}^{n}S_{1,\lambda}^{Y}(j,r)\frac{1}{j!}\big(\Delta^{j}q(1)+\Delta^{j}q(0) \big) \\
&=\frac{1}{2}\sum_{j=r}^{n}\sum_{i=0}^{j}(-1)^{j-i}\frac{1}{j!}\binom{j}{i}S_{1,\lambda}^{Y}(j,r)\big(q(i+1)+q(i)\big) \\
&=\frac{1}{2}\sum_{j=r}^{n}\sum_{l=j}^{n}S_{1,\lambda}^{Y}(j,r)S_{2}(l,j)\frac{1}{l!}\big(q^{(l)}(1)+q^{(l)}(0) \big).
\end{align*}
\end{theorem}
\begin{remark}
We have established the necessity of inverse relations for expressing arbitrary polynomials in terms of the probabilistic Euler polynomials, $\mathcal{E}_{n}^{Y}(x)$, and the probabilistic degenerate Euler polynomials, $\mathcal{E}_{n,\lambda}^{Y}(x)$. Theorems 2.1 through 2.7 provide explicit formulas for $S_{1}^{Y}(j,r)$ and $S_{1,\lambda}^{Y}(j,r)$ for several discrete and continuous random variables $Y$, thus enabling the explicit determination of these polynomial representations for those specific random variables. A more comprehensive analysis of such representations will be presented in subsequent publications.
\end{remark}

\section{\bf Conclusion}
The study of Stirling numbers (of both kinds) and their degenerate versions has been extended to a probabilistic context.  This involves defining probabilistic Stirling numbers of the first and second kind, both in standard and degenerate forms, associated with a random variable Y.  These probabilistic versions are constructed to maintain important properties like orthogonality and inverse relations, which are crucial for applications like expressing polynomials in terms of probabilistic Euler polynomials.  This builds upon prior work in defining probabilistic special numbers and polynomials, demonstrating a broader trend of incorporating probability into the study of these mathematical objects.

\vspace{0.3cm}

{\bf Competing interests.}
 The authors have no competing interests to declare that are relevant to the content of
this article and the authors contributed equally to this work.


\begin{thebibliography}{9}

\bibitem{1}
Adell, J. A. \emph{Probabilistic Stirling Numbers of the Second Kind and Applications,} J. Theor. Probab. 2022, 35: 636–652. https://doi.org/10.1007/s10959-020-01050-9
\bibitem{2}
Adell, J. A.; B{\'e}nyi, B. \emph{Probabilistic Stirling numbers and applications,} Aequat. Math. \textbf{98} (2024), 1627–1646.
\bibitem{3}
Adell, J. A.; Lekuona, A. \emph{A probabilistic generalization of the Stirling numbers of the second kind,} J. Number Theory \textbf{194} (2019), 225-355.
\bibitem{4}
Carlitz, L. \emph{Degenerate Stirling, Bernoulli and Eulerian numbers,} Utilitas Math. \textbf{15} (1979), 51-–88.
\bibitem{5}
Carlitz, L. \emph{A degenerate Staudt-Clausen theorem,} Arch. Math. (Basel) \textbf{7} (1956), 28-–33.
\bibitem{6}
Chen, L.; Kim, T.; Kim, D. S.; Lee, H.; Lee, S.-H. \emph{Probabilistic degenerate central Bell polynomials,} Math. Comput. Model. Dyn. Syst. \textbf{30} (2024), no. 1, 523-542. DOI: 10.1080/13873954.2024.2358899
\bibitem{7}
Comtet, L. \emph{Advanced combinatorics,} The art of finite and infinite expansions. Revised and enlarged edition. D. Reidel Publishing Co., Dordrecht, 1974.
\bibitem{8}
D. V. Dolgy, T. Kim, \emph{Some explicit formulas of degenerate Stirling numbers associated with the degenerate special numbers and polynomials,} Proc. Jangjeon Math. Soc. \textbf{21} (2018), no. 2, 309-317.
\bibitem{9}
Kim, D. S.; Kim, T. \emph{A note on new type degenerate Bernoulli numbers,} Russ. J.  Math. Phys. \textbf{27} (2020), no. 2, 227–235.
\bibitem{10}
Kim, D. S.; Kim, T. \emph{Degenerate Sheffer sequences and $\lambda$-Sheffer sequences,} J. Math. Anal. Appl. \textbf{493} (2021), 124521.
\bibitem{11}
Kim, D. S.; Kim, T. \emph{Stirling numbers associated with sequences of polynomials,} Appl. Comput. Math. \textbf{22} (2023), no. 1, 80-115.
\bibitem{12}
Kim, D. S.; Kim, T. \emph{Central factorial numbers associated with sequences of polynomials,} Math. Methods Appl. Sci. \textbf{46} (2023), no. 9, 10348-10383.
\bibitem{13}
Kim, D. S.; Kim, T.; Jang, G. W. \emph{A note on degenerate Stirling numbers of the first
kind,} Proc. Jangjeon Math. Soc. \textbf{21} (2018), 393–404.
\bibitem{14}
Kim, T. \emph{A note on degenerate Stirling polynomials of the second kind,} Proc. Jangjeon Math. Soc. \textbf{20}  (2017),  no. 3, 319-331.
\bibitem{15}
Kim, T.; Kim, D. S. \emph{Degenerate Laplace transform and degenerate gamma function,} Russ. J. Math. Phys. \textbf{24} (2017), no. 2, 241-248.
\bibitem{16}
Kim, T.; Kim, D. S. \emph{Note on the degenerate gamma function,} Russ. J. Math. Phys. \textbf{27} (2020), no. 3, 352–358.
\bibitem{17}
Kim, T.; Kim, D. S. \emph{Probabilistic degenerate Bell polynomials associated with random variables,} Russ. J. Math. Phys. \textbf{30} (2023), no. 4, 528–542.
\bibitem{18}
Kim, T.; Kim, D. S. \emph{Probabilistic degenerate Stirling numbers of the first kind and their applications,} Eur. J. Math. \textbf{10} (2024), Article No. 73.
\bibitem{19}
Kim, T.; Kim, D. S. \emph{Probabilistic Bernoulli and Euler polynomials,} Russ. J. Math.
Phys. \textbf{31} (2024), no. 1, 94-105.
\bibitem{20}
Kim, T.; Kim, D. S.; Kwon, H.-I. \emph{A note on degenerate Stirling numbers and their applications,} Proc. Jangjeon Math. Soc. \textbf{21} (2018), no. 2, 195-203.
\bibitem{21}
Kim, T.; Kim, D.S.; Kwon, J. \emph{Probabilistic degenerate Stirling polynomials of the second kind and their applications,} Math. Comput. Model. Dyn. Syst. \textbf{30} (2024), no. 1, 16-30.
\bibitem{22}
Kim, T.; Kim, D. S.; Lee, H.; Park, J.-W. \emph{A note on degenerate $r$-Stirling numbers,}      J. Inequal. Appl. (2020) 2020:225.
\bibitem{23}
Quaintance, J.; Gould, H. W. \emph{Combinatorial identities for Stirling numbers, The unpublished  notes of H. W. Gould, With a foreword by George E. Andrews,} World Scientific Publishing Co. Pte. Ltd., Singapore, 2016.
\bibitem{24}
Riordan, J. \emph{An introduction to combinatorial analysis,} Dover Publications, Inc., New York, 1958.
\bibitem{25}
Riordan, J. \emph{Combinatorial identities,} Wiley, New York, 1968.
\bibitem{26}
Roman, S. \emph{The umbral calculus,} Pure and Applied Mathematics, 111.Academic Press, Inc. [Harcourt Brace Jovanovich, Publishers], New York, 1984. x+193 pp. ISBN: 0-12-594380-6
\bibitem{27}
Ross, S. M. \emph{Introduction to probability models,} Twelfth edition, Academic Press, London, 2019.
\bibitem{28}
Ta, B.Q. \emph{Probabilistic Approach to Appell Polynomials,} Expo. Math. \textbf{33} (2015) no. 3, 269–294. https://doi.org/10.1016/j.exmath.2014.07.
\bibitem{29}
Xu, R.; Ma, Y.; Kim, T.; Kim, D. S.; Boulaaras, S. \emph{Probabilistic central Bell polynomials,} Fractals, to appear.
\end{thebibliography}
\end{document}